\documentclass[reqno]{amsproc}
\usepackage[latin1]{inputenc}
\usepackage{amssymb}
\usepackage{hyperref}
\usepackage{amsmath}
\usepackage{latexsym}
\usepackage{cite}
\usepackage{amssymb}
\usepackage{amsfonts}
\usepackage{amscd}
\usepackage{xcolor}
\usepackage{amstext,amsmath,amssymb,amsfonts}
\newtheorem{theorem}{Theorem}

\newtheorem{definition}[theorem]{Definition}
\newtheorem{example}[theorem]{Example}

\newtheorem{lemma}[theorem]{Lemma}

\newtheorem{proposition}[theorem]{Proposition}
\newtheorem{remark}[theorem]{Remark}

\newcommand{\A}{\mathcal{A}}

\newcommand{\B}{\mathcal{B}}

\newcommand{\M}{\mathcal{M}}

\newcommand{\beq}{\begin{eqnarray}}
\newcommand{\eeq}{\end{eqnarray}}
\newcommand{\beqs}{\begin{eqnarray*}}
\newcommand{\eeqs}{\end{eqnarray*}}
\newcommand{\bpro}{\begin{pro}}
\newcommand{\epro}{\end{pro}}
\newcommand{\blem}{\begin{lem}}
\newcommand{\elem}{\end{lem}}
\newcommand{\bdfn}{\begin{dfn}}
\newcommand{\edfn}{\end{dfn}}
\newcommand{\bcor}{\begin{cor}}
\newcommand{\ecor}{\end{cor}}
\newcommand{\bthm}{\begin{thm}}
\newcommand{\ethm}{\end{thm}}
\newcommand{\bex}{\begin{ex}}
\newcommand{\eex}{\end{ex}}
\newcommand{\brmk}{\begin{rmk}}
\newcommand{\ermk}{\end{rmk}}
\newcommand{\bpr}{\begin{pr}}
\newcommand{\epr}{\end{pr}}
\newcommand{\benum}{\begin{enumerate}}
\newcommand{\eenum}{\end{enumerate}}
\newcommand{\bitem}{\begin{itemize}}
\newcommand{\eitem}{\end{itemize}}

\chardef\bslash=`\\
\numberwithin{equation}{section}
\numberwithin{table}{section}
\numberwithin{theorem}{section}
\DeclareMathOperator{\id}{id}

\title[(Co)associative $3$-ary (co)algebras and   bialgebras]{(Co)associative $3$-ary (co)algebras and   infinitesimal bialgebras: construction and main properties \footnote{Preprint: ICMPA-MPA/2016/08 } }

\author{Mahouton Norbert Hounkonnou$^\ast$}
\address[$\ast$]{University of Abomey-Calavi,
International Chair in Mathematical Physics and Applications,
ICMPA-UNESCO Chair, 072 BP 50, Cotonou, Rep. of Benin}
\email{norbert.hounkonnou@cipma.uac.bj, with copy to hounkonnou@yahoo.fr}
\author{Gb\^ev\`ewou Damien  Houndedji$^\dagger$}
\address[$\dagger$]{University of Abomey-Calavi,
International Chair in Mathematical Physics and Applications,
ICMPA-UNESCO Chair, 072 BP 50, Cotonou, Rep. of Benin}
\email{ houndedjid@gmail.com}
\begin{document}
\maketitle

\today

\bigskip
\begin{abstract}
 The  (co)associative, partially (co)associative  and totally (co)associative $3$-ary (co) algebras and  infinitesimal bialgebras are constructed and discussed.  Their trimodules and matched pairs are defined and completely characterized. The main structural properties and relations are also deduced and analyzed.  
\end{abstract}

\section{Introduction}
An $n$-ary algebra is a  linear space endowed with an internal composition law involving $n$ elements:
 $\mu: V^{\otimes n}\rightarrow V.$ These  $n$-ary
   algebras, $n\ge 3,$ knew a versatile development since the discovery of the Nambu mechanics in 1973 (see \cite{[B22]}) and the work by S. Okubo  \cite{[B23]} on Yang-Baxter equation. The $n$-ary products were also defined by cubic matrices and a generalization of the notion of  determinant, called  hyperdeterminant, first introduced by Cayley in 1840, then rediscovered and generalized by Sokolov in 1972 \cite{sokolov},  and, still later, by Kapranov, Gelfand and Zelevinskii in 1994 \cite{[B12]}.

 The ternary algebraic structures  are particularly of a great  potential application in various domains of physics  and mathematics, and data processing.  Their subclass  known as Bagger-Lambert algebras \cite{[B3]} are involved in string theory and $M-$branes. 
A good compilation of their applications  can be found in the work by Kerner (see \cite{[B12]}-- \cite{[B15]}) on ternary
 and non-associative structures and their applications in physics. This author  investigated the use of
        $Z_{3}$-graded structures instead of $Z_{2}$-graded  in physics,  leading to 
         interesting results  in the construction of gauge theories.

 The $n$-ary algebras of associative type were studied by Lister, Loos, Myung and Carlsson (see
\cite{[B7]}, \cite{[B8]}, \cite{[B18]}, \cite{[B19]}). Relatively to their structure, these algebras encompass  two main classes: totally associative $n$-ary algebras and partially associative $n$-ary algebras, which also generate  other interesting variants of algebras. 
\begin{definition}
A totally associative $3$-ary algebra is  a $\mathcal{K}$-vector space $\mathcal{T}$ endowed with a trilinear operation $\mu$ satisfying, for all $x_{1}, x_{2}, x_{3}, x_{4}, x_{5}\in \mathcal{A}:$
\beq
\mu(\mu(x_{1}\otimes x_{2}\otimes x_{3})\otimes x_{4}\otimes x_{5})&=& \mu(x_{1}\otimes \mu(x_{2}\otimes x_{3}\otimes x_{4})\otimes x_{5})\cr &=& \mu(x_{1}\otimes x_{2}\otimes \mu(x_{3}\otimes x_{4}\otimes x_{5})).
\eeq 
\end{definition}
\begin{example}\label{et1}
Let $\mathcal{T}$ be a $2$-dimensional space vector with a basis $\lbrace e_{1}, e_{2} \rbrace$. The trilinear product $\mu$ on $\mathcal{T}$ defined by
\beqs
&&\mu(e_{1}\otimes e_{1}\otimes e_{1})= e_{1} \ \ \ \mu(e_{2}\otimes e_{2}\otimes e_{1})= e_{1} + e_{2}\cr
&&\mu(e_{1}\otimes e_{1}\otimes e_{2})= e_{2} \ \ \ \mu(e_{2}\otimes e_{2}\otimes e_{2})= e_{1} + 2e_{2}\cr
&&\mu(e_{1}\otimes e_{2}\otimes e_{1})= e_{2} \ \ \ \mu(e_{1}\otimes e_{2}\otimes e_{2})= e_{1} + e_{2}\cr
&&\mu(e_{2}\otimes e_{1}\otimes e_{1})= e_{2} \ \ \ \mu(e_{2}\otimes e_{1}\otimes e_{2})= e_{1} + e_{2}
\eeqs
defines a totally associative $3$-ary algebra.
\end{example}

\begin{definition}
A weak totally associative $3$-ary algebra is  a $\mathcal{K}$-vector space $\mathcal{W}$ equipped with a trilinear operation $\mu$ satisfying, for all $x_{1}, x_{2}, x_{3}, x_{4}, x_{5}\in \mathcal{W}:$
\beq
\mu(\mu(x_{1}\otimes x_{2}\otimes x_{3})\otimes x_{4}\otimes x_{5})= \mu(x_{1}\otimes x_{2}\otimes \mu(x_{3}\otimes x_{4}\otimes x_{5})).
\eeq
\end{definition}
\begin{remark}
Naturally, any totally associative $3$-ary algebra is a weak totally associative $3$-ary algebra.
\end{remark}
\begin{definition}
A partially associative $3$-ary algebra is  a $\mathcal{K}$-vector space $\mathcal{P}$ with a trilinear operation $\mu$ satisfying, for all $x_{1}, x_{2}, x_{3}, x_{4}, x_{5}\in \mathcal{P}:$
\beq \label{qp1}
\mu(\mu(x_{1}\otimes x_{2}\otimes x_{3})\otimes x_{4}\otimes x_{5})+ \mu(x_{1}\otimes \mu(x_{2}\otimes x_{3}\otimes x_{4})\otimes x_{5})+\cr
 \mu(x_{1}\otimes x_{2}\otimes \mu(x_{3}\otimes x_{4}\otimes x_{5}))=0.
\eeq 
\end{definition}
\begin{example}\label{ep1}
Let $\mathcal{P}$ be a $2$-dimensional vector space  with a basis $\lbrace e_{1}, e_{2} \rbrace$. The trilinear product $\mu$ on $\mathcal{P}$ determined by
\beqs
\mu(e_{1}\otimes e_{1}\otimes e_{1})= e_{2}; \  \mu(e_{i}\otimes e_{j}\otimes e_{k})=0,\ i, j, k=1, 2 \mbox{ with } (i, j, k)\neq (1, 1, 1).
\eeqs
defines a partially associatve $3$- ary algebra.
\end{example}
\begin{remark} \label{r1}
Let $(\mathcal{A}, m)$ be a bilinear associative algebra. Then, the trilinear operation defined on $\mathcal{A}$ by 
\beqs
\mu(x_{1}\otimes x_{2}\otimes x_{3})= m(m(x_{1}\otimes x_{2})\otimes x_{3})
\eeqs
produces on the vector space $\mathcal{A}$ a structure of totally associative $3$-ary algebra which is not partially associative.
\end{remark}
From the Remark \ref{r1}, we can deduce a series of  examples of totally associative $3$-ary algebras, and hence of weak totally associative $3$-ary algebras. In the next example, we provide a class of  compatible $2-$dimensional  totally associative $3$-ary algebras derived from  $2-$dimensional associative algebras.
\begin{example}
Let $(\mathcal{A}, m)$ be an associative algebra with a basis $\lbrace e_{1}, e_{2}
   \rbrace.$ We denote by  $(\mathcal{A_i}, \mu), \,i=1,7,$ the related compatible totally associative $3$-ary algebras defined as follows:
\begin{itemize}
\item $\A_{1}$: $m(e_{1}\otimes e_{1})= e_{1}, m(e_{2}\otimes e_{2})= e_{2};$\ \  $\mu(e_{1}\otimes e_{1}\otimes e_{1})= e_{1}, \mu(e_{2}\otimes e_{2}\otimes e_{2})= e_{2}$.
\item $\A_{2}$: $m(e_{2}\otimes e_{2})= e_{2},$ $m(e_{1}\otimes e_{2})= m(e_{2}\otimes e_{1})= e_{1};$ \ \  $\mu(e_{2}\otimes e_{2}\otimes e_{2})= e_{2},$ \\
 $\mu(e_{1}\otimes e_{2}\otimes e_{2})=$ $\mu(e_{2}\otimes e_{1}\otimes e_{2})= \mu(e_{2}\otimes e_{2}\otimes e_{1})= e_{1}.$
\item $\A_{3}$: $m(e_{1}\otimes e_{1})= e_{1};$ \ \  $\mu(e_{1}\otimes e_{1}\otimes e_{1})= e_{1}.$
\item $\A_{4}$: $m(e_{i}\otimes e_{j})= 0, i, j=1, 2;$ \ \ $\mu(e_{i}\otimes e_{j}\otimes e_{k})= 0, i, j, k=1, 2.$
\item $\A_{5}$: $m(e_{1}\otimes e_{1})= e_{2};$ \ \ $\mu(e_{i}\otimes e_{j}\otimes e_{k})= 0, i, j, k=1, 2.$
\item $\A_{6}$:  $m(e_{2}\otimes e_{1})= e_{1}, m(e_{2}\otimes e_{2})= e_{2};$ \ \ $\mu(e_{2}\otimes e_{2}\otimes e_{1})= e_{1}, \mu(e_{2}\otimes e_{2}\otimes e_{2})= e_{2}.$
\item $\A_{7}$: $m(e_{1}\otimes e_{2})= e_{1}, m(e_{2}\otimes e_{2})= e_{2};$ \ \ $\mu(e_{1}\otimes e_{2}\otimes e_{2})= e_{1}, \mu(e_{2}\otimes e_{2}\otimes e_{2})= e_{2}.$
\end{itemize}
\end{example}

Note that  the $3$-ary algebras  given by subspaces of an
 associative algebra, closed  under the ternary product $(x, y, z)\mapsto xyz,$
 are linked to the ternary operation defined by Hestenes 
 \cite{[Hestenes]}  on a linear space of rectangular matrices $A, B, C\in \M_{m, n},$ with
 complex entries by $AB^{\ast}C,$ where $B^{\ast}$ is the conjugate transpose matrix of $B.$
 This operation, strictly speaking, does not define a ternary algebra product on $\M_{m, n}$ as it is
 linear on the first and the third arguments, but conjugate-linear on the second argument.
  It satisfies identities, sometimes referred to as identities of total associativity of
 second kind, which only slightly differ from the identities of totally associative
 algebras. 

The totally associative ternary algebras are also sometimes called
 associative triple systems. 

The cohomology of totally associative $n$-ary algebras was
 studied by Carlsson through the embedding  \cite{[B8]}. In \cite{[B1]},  the $1$-parameter formal deformation theory was extended to ternary algebras of associative  type, while  in \cite{[B2]} discussions were made on their 
 cohomologies in the context of  deformations.  See also \cite{[B11]}, \cite{[B9]}, and  \cite{[B25]} (and references therein ).

The extension of the notion of associativity to $n-$ary product is not trivial. The most natural procedure might be based on  the notion of totally associativity. Unfortunately, this notion is not auto-dual in the
operadic point of view. So,  it is necessary to introduce a most general notion of
associativity under the concept of   partial associativity.   

In this paper, we focus on  the  case $n=3.$  Of course, the structure of classes of
associative $n$-ary algebras,(for $n\geq 3$), is more complicated than that  of associative algebras. Hence,  their exhaustive investigation in order to derive their relevant relationships is of some importance in algebra. 
Besides, the concepts of associative coalgebras and bialgebras are fundamental  in the theory of 
associative algebras. In \cite{[M1]} and
 \cite{[M2]},  Aguiar  developed the basic theory of infinitesimal bialgebrass and infinitesimal
  hopf bialgebras. An infinitesimal bialgebra is at the same time an algebra and a coalgebra, in such
   a way that the comultiplication is a derivation. Aguiar established many properties of ordinary Hopf algebras which possess infinitesimal version. He introduced bicrossproducts, quasitriangular
    infinitesimal bialgebras, the corresponding infinitesimal Yang-Baxter equation and a notion of Drinfeld's double for infinitesimal Hopf algebras. He also  showed that non degenerate antisymmetric solutions of associative Yang-Baxter equations are in one-to-one correspondence with non degenerate cyclic $2$-cocycles. Furthermore,  Bai  established a clear analogy between the antisymmetric infinitesimal bialgebras and the dendriform $D$-bialgebras \cite{[C.Bai]}. Motivated by all these studies, we define, in this work,  the concepts of  totally and partially coassociative $3$-ary coalgebras and  totally and partially associative $3$-ary  infinitesimal bialgebras, and investigate their main properties as well as their relationships with associative $3$-ary infinitesimal bialgebras.

The paper is organized as follows. In section $2$, we construct partially and totally coassociative $3$-ary coalgebras and discuss their main properties. Section $3$ is devoted to the construction of trimodules and matched pairs of totally and partially associative $3-$ary algebras. Then, in Section $4$, we define the partially and totally associative $3$-ary infinitesimal bialgebras and investigate  their relation  with associative $3$-ary algebras. Section $5$ is devoted to concluding remarks.

%

 Throughout the  paper, we consider a fixed field $\mathcal{K}$ of characteristic zero.
\section{Coassociative $3$-ary coalgebras}
In this section, we introduce and develop the concepts of totally coassociative $3$-ary coalgebra, weak totally coassociative $3$-ary coalgebra and partially coassociative $3$-ary coalgebra. 
\subsection{Definitions}
Let us start with the following definitions.
\begin{definition}
A totally associative $3$-ary algebra is  a $\mathcal{K}$-vector space $\mathcal{T}$ equipped with a trilinear operation $\mu$ satisfying
\beq
\mu\circ(\mu\otimes \id\otimes \id)= \mu\circ(\id\otimes \mu\otimes \id)= \mu\circ(\id\otimes \id\otimes \mu).
\eeq 
\end{definition}
 \begin{definition}
A weak totally associative $3$-ary algebra is a $\mathcal{K}$-vector space $\mathcal{W}$ with a  trilinear operation $\mu$ satisfying
\beq
\mu\circ(\mu\otimes \id\otimes \id)= \mu\circ(\id\otimes \id\otimes \mu).
\eeq
\end{definition}
\begin{definition}
A partially associative $3$-ary algebra is  a $\mathcal{K}$-vector space $\mathcal{P}$ endowed with a trilinear operation $\mu$ satisfying
\beq
\mu\circ(\mu\otimes \id\otimes \id +\id\otimes \mu\otimes \id +\id\otimes \id\otimes \mu)=0.
\eeq 
\end{definition}
Let $(\mathcal{A}, \mu)$ be an associative  $3$-ary algebra and $\mathcal{A}^{\ast}$ be its dual space. Then, we get the dual mapping $\mu^{\ast}: \mathcal{A}^{\ast}\rightarrow \mathcal{A}^{\ast}\otimes\mathcal{A}^{\ast}\otimes\mathcal{A}^{\ast}$ of $\mu$, for every $x, y, z\in \mathcal{A}$ and $\xi, \eta, \gamma \in \mathcal{A}^{\ast}$, 
\beq \label{qdual}
\langle \mu^{\ast}(\xi), x\otimes y\otimes z\rangle= \langle\xi, \mu(x, y, z) \rangle,
\eeq 
\beq
\langle \xi\otimes\eta\otimes\gamma, x\otimes y\otimes z\rangle= \langle\xi, x\rangle \langle\eta, y\rangle \langle\gamma, z\rangle,
\eeq
where $\langle, \rangle$ is the natural nondegenerate symmetric bilinear form on the vector space $\mathcal{A}\oplus \mathcal{A}^{\ast}$ defined by $\langle \xi, x\rangle= \xi(x), \xi\in \mathcal{A}^{\ast}, x\in \mathcal{A}$.

Further, by the definition of associative $3$-ary algebras, we have $Im(\mu^{\ast})\subseteq \mathcal{A}^{\ast}\otimes\mathcal{A}^{\ast}\otimes\mathcal{A}^{\ast},$ and
\beqs
(\mu^{\ast}\otimes \id\otimes \id)\circ\mu^{\ast}= (\id\otimes \mu^{\ast}\otimes \id)\circ\mu^{\ast}= (\id\otimes \id\otimes \mu^{\ast})\circ\mu^{\ast} 
\eeqs 
and
\beqs
(\mu^{\ast}\otimes \id\otimes \id +\id\otimes \mu^{\ast}\otimes \id+ \id\otimes \id\otimes \mu^{\ast})\circ\mu^{\ast}=0, 
\eeqs 
for the totally associative $3$-ary algebras and the partially associative $3$-ary algebras, respectively. That is, for every $x_{1}, x_{2}, x_{3}, x_{4}, x_{5}\in \mathcal{A}$ and $\xi\in \mathcal{A}^{\ast}$, we have
\beqs
&&\langle(\mu^{\ast}\otimes\id\otimes\id)\circ\mu^{\ast}(\xi), x_{1}\otimes x_{2}\otimes x_{3}\otimes x_{4}\otimes x_{5}\rangle\cr
&&=\langle\xi, \mu\circ(\mu\otimes\id\otimes\id)(x_{1}\otimes x_{2}\otimes x_{3}\otimes x_{4}\otimes x_{5})\rangle\cr
&&=\langle\xi, \mu(\mu(x_{1}\otimes x_{2}\otimes x_{3})\otimes x_{4}\otimes x_{5})\rangle \cr
&&=\langle\xi, \mu(x_{1}\otimes \mu(x_{2}\otimes x_{3}\otimes x_{4})\otimes x_{5})\rangle\cr
&&=\langle(\id\otimes\mu^{\ast}\otimes\id)\circ\mu^{\ast}(\xi), x_{1}\otimes x_{2}\otimes x_{3}\otimes x_{4}\otimes x_{5}\rangle\cr
&&=\langle\xi, \mu(x_{1}\otimes x_{2}\otimes \mu(x_{3}\otimes x_{4}\otimes x_{5}))\rangle\cr
&&=\langle(\id\otimes\id\otimes\mu^{\ast})\circ\mu^{\ast}(\xi), x_{1}\otimes x_{2}\otimes x_{3}\otimes x_{4}\otimes x_{5}\rangle\cr 
\eeqs 
and
\beqs
&&\langle(\mu^{\ast}\otimes \id\otimes \id +\id\otimes \mu^{\ast}\otimes \id+ \id\otimes \id\otimes \mu^{\ast})\circ\mu^{\ast}(\xi), x_{1}\otimes x_{2}\otimes x_{3}\otimes x_{4}\otimes x_{5}\rangle\cr
=&&\langle\xi, \mu\circ(\mu\otimes \id\otimes \id +\id\otimes \mu\otimes \id +\id\otimes \id\otimes \mu)(x_{1}\otimes x_{2}\otimes x_{3}\otimes x_{4}\otimes x_{5})\rangle,
\eeqs
respectively.

Provided the above, the following definitions are in order.
 \begin{definition}\label{Def1}
 A totally coassociative $3$-ary coalgebra $(\mathcal{T}, \Delta)$ is a vector space $\mathcal{T}$ with a linear mapping $\Delta: \mathcal{T}\rightarrow \mathcal{T}\otimes\mathcal{T}\otimes\mathcal{T}$ satisfying
 \beq \label{totally coassociativity}
 (\Delta\otimes \id\otimes \id)\circ\Delta= (\id\otimes \Delta\otimes \id)\circ\Delta= (\id\otimes \id\otimes \Delta)\circ\Delta. 
 \eeq
  \end{definition}
The Eq. \ref{totally coassociativity} is also called standard ternary coassociativity condition  in \cite{[S.Duplij]}.

 \begin{definition}\label{Def2}
 A weak totally coassociative $3$-ary coalgebra $(\mathcal{W}, \Delta)$ is a vector space $\mathcal{W}$ with a linear mapping $\Delta: \mathcal{W}\rightarrow \mathcal{W}\otimes\mathcal{W}\otimes\mathcal{W}$ satisfying
 \beq 
 (\Delta\otimes \id\otimes \id)\circ\Delta= (\id\otimes \id\otimes \Delta)\circ\Delta.
 \eeq
  \end{definition}  
  
 \begin{definition}\label{Def3}
 A partially coassociative $3$-ary coalgebra $(\mathcal{P}, \Delta)$ is a vector space $\mathcal{P}$ with a linear mapping $\Delta: \mathcal{P}\rightarrow \mathcal{P}\otimes\mathcal{P}\otimes\mathcal{P}$ satisfying
 \beq 
 (\Delta\otimes \id\otimes \id +\id\otimes \Delta\otimes \id +\id\otimes \id\otimes \Delta)\circ\Delta= 0. 
 \eeq
\end{definition}
\subsection{Main results}

Let us now  examine the  partially coassociative $3$-ary coalgebras in terms of structure constants. For that,  let $(\mathcal{A}, \Delta)$ be a partially coassociative $3$-ary coalgebra with a basis $e_{1}, ....., e_{n}$. Assume that 
\beq\label{q1}
\Delta(e_{l})= \sum_{1\leq r, s, t\leq n}c^{l}_{rst}e_{r}\otimes e_{s}\otimes e_{t}, \ \ c^{l}_{rst}\in \mathcal{K}, 1\leq l\leq n.
\eeq
Then, we obtain
\beqs
&&(\Delta\otimes\id\otimes\id + \id\otimes\Delta\otimes\id + \id\otimes\id\otimes\Delta)\circ\Delta(e_{l})\cr
&&=(\Delta\otimes\id\otimes\id + \id\otimes\Delta\otimes\id + \id\otimes\id\otimes\Delta)\left( \sum_{1\leq r, s, t\leq n}c^{l}_{rst}e_{r}\otimes e_{s}\otimes e_{t}\right) \cr
&&=\sum_{1\leq r, s, t\leq n}c^{l}_{rst}\left( \Delta(e_{r})\otimes e_{s}\otimes e_{t} + e_{r}\otimes \Delta(e_{s})\otimes e_{t} + e_{r}\otimes e_{s}\otimes \Delta(e_{t})\right) \cr
&&=\sum_{1\leq r, s, t\leq n}\sum_{1\leq i, j, k\leq n}c^{l}_{rst}c^{r}_{ijk}e_{i}\otimes e_{j}\otimes e_{k}\otimes e_{s}\otimes e_{t} +\cr
&&\sum_{1\leq r, s, t\leq n}\sum_{1\leq i, j, k\leq n}c^{l}_{rst}c^{s}_{ijk}e_{r}\otimes e_{i}\otimes e_{j}\otimes e_{k}\otimes e_{t} + \cr 
&&\sum_{1\leq r, s, t\leq n}\sum_{1\leq i, j, k\leq n}c^{l}_{rst}c^{t}_{ijk}e_{r}\otimes e_{s}\otimes e_{i}\otimes e_{j}\otimes e_{k}\cr
&&=\sum_{1\leq r, s, t\leq n}\sum_{1\leq i, j, k\leq n}(c^{l}_{rst}c^{r}_{ijk}+ c^{l}_{irt}c^{r}_{jks} + c^{l}_{ijr}c^{r}_{kst})e_{i}\otimes e_{j}\otimes e_{k}\otimes e_{s}\otimes e_{t}
\eeqs 
leading to
\beq \label{q2}
\sum^{n}_{r=1}(c^{l}_{rst}c^{r}_{ijk}+ c^{l}_{irt}c^{r}_{jks} + c^{l}_{ijr}c^{r}_{kst})=0, \ \ 1\leq i, j, k, s, t, l\leq n.
\eeq
By a similar discussion, for a totally coassociative $3$-ary coalgebra, we get
\beq \label{q3}
\sum^{n}_{r=1}c^{l}_{rst}c^{r}_{ijk}= \sum^{n}_{r=1}c^{l}_{irt}c^{r}_{jks}=\sum^{n}_{r=1} c^{l}_{ijr}c^{r}_{kst}, \ \ 1\leq i, j, k, s, t, l\leq n,
\eeq
while for a weak totally coassociative $3$-ary coalgebra, 
\beq \label{q4}
\sum^{n}_{r=1}c^{l}_{rst}c^{r}_{ijk}= \sum^{n}_{r=1} c^{l}_{ijr}c^{r}_{kst}, \ \ 1\leq i, j, k, s, t, l\leq n.
\eeq
Therefore, we infer the following statement.
 \begin{theorem}\label{theo1}
 Let $\mathcal{A}$ be an $n$-dimensional vector space with a basis $e_{1},....., e_{n},$ and
 
  $\Delta: \mathcal{A} \rightarrow \mathcal{A}\otimes \mathcal{A}\otimes \mathcal{A}$ be defined as (\ref{q1}). Then,
 \begin{enumerate}
 \item $(\mathcal{A}, \Delta)$ is a partially coassociative $3$-ary coalgebra if and only if the constants $c^{l}_{ijk},$\\
 $1\leq i, j, k\leq n$ satisfy the identity (\ref{q2});
 \item $(\mathcal{A}, \Delta)$ is a totally coassociative $3$-ary coalgebra if and only if the constants $c^{l}_{ijk},$ \\
 $1\leq i, j, k\leq n$ satisfy the identity (\ref{q3});
 \item $(\mathcal{A}, \Delta)$ is a weak totally coassociative $3$-ary coalgebra if and only if the constants $c^{l}_{ijk},$ \\
 $1\leq i, j, k\leq n$ satisfy the identity (\ref{q4}).
 \end{enumerate}
  \end{theorem}
  
  Now, let $(\mathcal{A}, \mu)$ be a partially associative $3$-ary algebra with a basis $e_{1}, e_{2},....., e_{n},$ and the mutiplication $\mu$ of $\mathcal{A}$ in this basis be defined as follows:
  
\beq
\mu(e_{r}, e_{s}, e_{t})=\sum^{n}_{l=1}c^{l}_{rst}e_{l}, \ \  c^{l}_{rst}\in \mathcal{K}, \ \ 1\leq r, s, t\leq n.
\eeq  
Using the condition (\ref{qp1}), we have:
\beqs
&&\mu(\mu(e_{r}, e_{s}, e_{t}), e_{i}, e_{j}) + \mu(e_{r}, \mu(e_{s}, e_{t}, e_{i}), e_{j}) + \mu(e_{r}, e_{s}, \mu(e_{t}, e_{i}, e_{j}))\cr
&&=\mu( \sum^{n}_{l=1}c^{l}_{rst} e_{l}, e_{i}, e_{j}) + \mu(e_{r}, \sum^{n}_{l=1}c^{l}_{sti}e_{l}, e_{j}) + \mu(e_{r}, e_{s}, \sum^{n}_{l=1}c^{l}_{tij}e_{l})\cr
&&= \sum^{n}_{k=1}\sum^{n}_{l=1}c^{l}_{rst}c^{k}_{lij} e_{k} + \sum^{n}_{k=1}\sum^{n}_{l=1}c^{l}_{sti}c^{k}_{rlj} e_{k} + \sum^{n}_{k=1}\sum^{n}_{l=1}c^{l}_{tij}c^{k}_{rsl} e_{k}\cr
&&=0
\eeqs
 yielding $\sum^{n}_{l=1}(c^{l}_{rst}c^{k}_{lij} + c^{l}_{sti}c^{k}_{rlj} +
   c^{l}_{tij}c^{k}_{rsl})=0,$ i. e.  $ \lbrace c^{l}_{i_{1}i_{2}i_{3}, 1\leq i_{1}, i_{2},
    i_{3}\leq n} \rbrace$ satisfies the identity (\ref{q2}).
    
    Similarly,  for a totally associative $3$-ary algebra and a weak totally associative $3$-ary algebra, we derive, respectively:
        
         $\sum^{n}_{l=1}c^{l}_{rst}c^{k}_{lij}= \sum^{n}_{l=1}c^{l}_{sti}c^{k}_{rlj}= 
   \sum^{n}_{l=1}c^{l}_{tij}c^{k}_{rsl},$ i. e. $ \lbrace c^{l}_{i_{1}i_{2}i_{3}, 1\leq i_{1}, i_{2}, i_{3}\leq n} \rbrace$ satisfies the identity (\ref{q3}),
   
   and
   
   $\sum^{n}_{l=1}c^{l}_{rst}c^{k}_{lij}= \sum^{n}_{l=1}c^{l}_{tij}c^{k}_{rsl},$ i. e.  $ \lbrace c^{l}_{i_{1}, i_{2}, i_{3}, 1\leq i_{1}, i_{2}, i_{3}\leq n} \rbrace$ satisfies the identity (\ref{q4}).
   
   Let $\mathcal{A}^{\ast}$ be the dual space of partially associative $3$-ary algebra
    $(\mathcal{A}, \mu)$, and $e^{\ast}_{1},....., e^{\ast}_{n}$ be the dual basis of $e_{1},....., e_{n}, \langle e^{\ast}_{1}, e_{j}\rangle= \delta_{ij}, 1\leq i, j\leq n.$ Assume that $\mu^{\ast}: \mathcal{A}^{\ast}$ $ \rightarrow$ $ \mathcal{A}^{\ast}\otimes \mathcal{A}^{\ast}\otimes \mathcal{A}^{\ast}$ is the dual mapping of $\mu$ defined by (\ref{qdual}). Then, for every $1\leq l\leq n,$ we have 
    \beq
    \mu^{\ast}(e^{\ast}_{l})= \sum_{1\leq r, s, t\leq n}c^{l}_{rst}e^{\ast}_{r}\otimes e^{\ast}_{s}\otimes e^{\ast}_{t}.
    \eeq
    Following the identities (\ref{q1}) and (\ref{q2}), ($\mathcal{A}^{\ast}, \mu^{\ast}$) is a partially coassociative $3$-ary coalgebra.  

Conversely, if $(\mathcal{A}, \Delta)$ is a partially coassociative $3$-ary coalgebra with a basis $e_{1}, ....., e_{n}$ satisfying (\ref{q1}), $\mathcal{A}^{\ast}$ is the dual space of $\mathcal{A}$ with the dual basis $ e^{\ast}_{1}, ...., e^{\ast}_{n}.$ Then, the dual mapping\\ $\Delta^{\ast}:$ $\mathcal{A}^{\ast}$ $ \rightarrow$ $ \mathcal{A}^{\ast}\otimes \mathcal{A}^{\ast}\otimes \mathcal{A}^{\ast}$ of $\Delta$ satisfies, for every $\xi, \eta, \gamma \in \mathcal{A}^{\ast}, x\in \mathcal{A},$
\beq \label{eqco}
\langle \Delta^{\ast}(\xi, \eta, \gamma), x\rangle= \langle\xi\otimes\eta\otimes\gamma, \Delta(x) \rangle.
\eeq 
Then, $\Delta^{\ast}(e^{\ast}_{r}, e^{\ast}_{s}, e^{\ast}_{t})= \sum^{n}_{l=1}c^{l}_{rst}e^{\ast}_{l}, c^{l}_{rst}\in \mathcal{K}, 1\leq r, s, t, l\leq n$ and $\Delta^{\ast}$ satisfies identity (\ref{q2}).
\begin{remark}
The above constructions and discussions on the totally associative $3$-ary algebra and  weak totally associative $3$-ary algebra also remain valid in this case. 
\end{remark}

Therefore,  the following results are true.
\begin{theorem} \label{teo-co-asso}
Let $\mathcal{A}$ be a vector space over a field $\mathcal{K}$, and $\Delta: \mathcal{A}\rightarrow \mathcal{A}\otimes \mathcal{A}\otimes \mathcal{A}.$ Then,
\begin{enumerate}
\item $(\mathcal{A}, \Delta)$ is a partially coassociative $3$-ary coalgebra if and only if $(\mathcal{A}^{\ast}, \Delta^{\ast})$ is a partially associative $3$-ary algebra.
\item $(\mathcal{A}, \Delta)$ is a totally coassociative $3$-ary coalgebra if and only if $(\mathcal{A}^{\ast}, \Delta^{\ast})$ is a totally associative $3$-ary algebra.
\item $(\mathcal{A}, \Delta)$ is a weak totally coassociative $3$-ary coalgebra if and only if $(\mathcal{A}^{\ast}, \Delta^{\ast})$ is a weak totally associative $3$-ary algebra.
\end{enumerate}
\end{theorem}

We can also give an equivalence description of (\ref{teo-co-asso}) as below.
\begin{theorem} \label{teo-co-asso1}
Let $\mathcal{A}$ be a vector space over a field $\mathcal{K}$, and $\mu: \mathcal{A}\otimes \mathcal{A}\otimes \mathcal{A} \rightarrow \mathcal{A}$ be a trilinear mapping. Then,
\begin{enumerate}
\item $(\mathcal{A}, \mu)$ is a partially associative $3$-ary algebra if and only if $(\mathcal{A}^{\ast}, \mu^{\ast})$ is a partially coassociative $3$-ary coalgebra.
\item $(\mathcal{A}, \mu)$ is a totally associative $3$-ary algebra if and only if $(\mathcal{A}^{\ast}, \mu^{\ast})$ is a totally coassociative $3$-ary coalgebra.
\item $(\mathcal{A}, \mu)$ is a weak totally associative $3$-ary algebra if and only if $(\mathcal{A}^{\ast}, \mu^{\ast})$ is a weak totally coassociative $3$-ary coalgebra.
\end{enumerate}
\end{theorem}

\begin{example}
Let $(\mathcal{T}^{\ast}, \mu^{\ast})$ be the dual of totally associative $3$-ary algebra $(\mathcal{T}, \mu)$ in\\ Example (\ref{et1}). The product $\mu^{\ast}$ on $\mathcal{T}^{\ast}$ is given by
\beqs
\mu^{\ast}(e^{\ast}_{1})&=& e^{\ast}_{1}\otimes e^{\ast}_{1}\otimes e^{\ast}_{1} + e^{\ast}_{1}\otimes e^{\ast}_{2}\otimes e^{\ast}_{2} + e^{\ast}_{2}\otimes e^{\ast}_{2}\otimes e^{\ast}_{1} + e^{\ast}_{2}\otimes e^{\ast}_{2}\otimes e^{\ast}_{2} + e^{\ast}_{2}\otimes e^{\ast}_{1}\otimes e^{\ast}_{2}\cr
\mu^{\ast}(e^{\ast}_{2})&=& e^{\ast}_{1}\otimes e^{\ast}_{1}\otimes e^{\ast}_{2} + e^{\ast}_{1}\otimes e^{\ast}_{2}\otimes e^{\ast}_{2} + e^{\ast}_{2}\otimes e^{\ast}_{1}\otimes e^{\ast}_{1} + e^{\ast}_{2}\otimes e^{\ast}_{2}\otimes e^{\ast}_{1} + e^{\ast}_{2}\otimes e^{\ast}_{1}\otimes e^{\ast}_{2} +\cr
&& e^{\ast}_{1}\otimes e^{\ast}_{2}\otimes e^{\ast}_{1} + 2e^{\ast}_{2}\otimes e^{\ast}_{2}\otimes e^{\ast}_{2}.
\eeqs
$(\mathcal{T}^{\ast}, \mu^{\ast})$ is a totally coassociative $3$-ary coalgebra.
\end{example}
\begin{example}
Let $(\mathcal{P}^{\ast}, \mu^{\ast})$ be the dual of partially associative $3$-ary algebra $(\mathcal{P}, \mu)$ in\\ Example (\ref{ep1}). The product $\mu^{\ast}$ on $\mathcal{P}^{\ast}$ is given by
\beqs
\mu^{\ast}(e^{\ast}_{2})= e^{\ast}_{1}\otimes e^{\ast}_{1}\otimes e^{\ast}_{1}; \   \mu(e^{\ast}_{1})=0.
\eeqs
$(\mathcal{P}^{\ast}, \mu^{\ast})$ is a partially coassociative $3$-ary coalgebra.
\end{example}
Let us recall that a coassossiative $3$-ary coalgebra is a partially coassossiative $3$-ary coalgebra, or a totally coassossiative $3$-ary coalgebra, or a weak totally coassossiative $3$-ary coalgebra.
\begin{definition}
Let ($\mathcal{A}_{1}, \Delta_{1}$) and ($\mathcal{A}_{2}, \Delta_{2}$) be two coassociative $3$-ary coalgebras. If there is a linear isomorphism $\varphi: \mathcal{A}_{1}\rightarrow \mathcal{A}_{2}$ satisfying
\beq
(\varphi\otimes\varphi\otimes\varphi)(\Delta_{1}(e))= \Delta_{2}(\varphi(e)), \mbox{ for every } e \in \mathcal{A}_{1},
\eeq 
then $(\mathcal{A}_{1}, \Delta_{1})$ is isomorphic to $(\mathcal{A}_{2}, \Delta_{2}),$ and $\varphi$ is called a coassociative $3$-ary coalgebra isomorphism, where 
\beq
(\varphi\otimes\varphi\otimes\varphi)\sum_{i}(a_{i}\otimes b_{i}\otimes c_{i})= \sum_{i}\varphi(a_{i})\otimes\varphi(b_{i})\otimes\varphi(c_{i}).
\eeq
\end{definition}
\begin{theorem}\label{theo2}
Let ($\mathcal{A}_{1}, \Delta_{1}$) and ($\mathcal{A}_{2}, \Delta_{2}$) be two coassociative $3$-ary
 coalgebras. Then, $\varphi: \mathcal{A}_{1}\rightarrow \mathcal{A}_{2}$ is a coassociative $3$-ary
  coalgebra isomorphism from ($\mathcal{A}_{1}, \Delta_{1}$) to ($\mathcal{A}_{2}, \Delta_{2}$) if and
   only if the dual mapping $ \varphi^{\ast}: \mathcal{A}^{\ast}_{2}\rightarrow
    \mathcal{A}^{\ast}_{1} $ is an associative $3$-ary algebra isomorphism from
     ($\mathcal{A}^{\ast}_{2}, \Delta^{\ast}_{2}$) to ($\mathcal{A}^{\ast}_{1}, \Delta^{\ast}_{1}$),
      where for every $\xi \in \mathcal{A}^{\ast}_{2}, v\in \mathcal{A}_{1}, \langle  \varphi^{\ast
      }(\xi), v\rangle= \langle\xi, \varphi(v)\rangle$.
\end{theorem}
\textbf{Proof:}

Let ($\mathcal{A}_{1}, \Delta_{1}$) and ($\mathcal{A}_{2}, \Delta_{2}$) be two coassociative $3$-ary
 coalgebras. It follows that  ($\mathcal{A}^{\ast}_{1}, \Delta^{\ast}_{1}$) and ($\mathcal{A}^{\ast}_{2},
 \Delta^{\ast}_{2}$) are two associative $3$-ary algebras. Let $\varphi: \mathcal{A}_{1}\rightarrow 
\mathcal{A}_{2}$ be a coassociative $3$-ary coalgebra isomorphism from ($\mathcal{A}_{1},
\Delta_{1}$) to ($\mathcal{A}_{2}, \Delta_{2}$). Hence, the dual mapping $\varphi^{\ast}:
 \mathcal{A}^{\ast}_{2}\rightarrow \mathcal{A}^{\ast}_{1}$ is a linear isomorphism and for every
  $\xi, \eta, \gamma \in \mathcal{A}^{\ast}_{2}$, $x \in \mathcal{A}_{1}:$
  \beqs
\langle \varphi^{\ast}\Delta^{\ast}_{2}(\xi, \eta, \gamma), x\rangle &=&\langle \xi\otimes\eta\otimes\gamma, \Delta_{2}(\varphi(x))\rangle\cr &=& \langle \xi\otimes\eta\otimes\gamma, (\varphi\otimes \varphi\otimes\varphi)\Delta_{1}(x)\rangle\cr
&=&\langle \varphi^{\ast}(\xi)\otimes\varphi^{\ast}(\eta)\otimes\varphi^{\ast}(\gamma), \Delta_{1}(x)\rangle \cr &=& \langle \Delta^{\ast}_{1}(\varphi^{\ast}(\xi), \varphi^{\ast}(\eta), \varphi^{\ast}(\gamma)), x\rangle.
  \eeqs
Then, $\varphi^{\ast}\Delta^{\ast}_{2}(\xi, \eta, \gamma)= \Delta^{\ast}_{1}(\varphi^{\ast}(\xi), \varphi^{\ast}(\eta), \varphi^{\ast}(\gamma)),$ that is, $\varphi^{\ast}$ is an associative $3$-ary algebra isomorphism.
$ \hfill \square $
\section{Trimodules and matched pairs of associative $3$-ary  algebras}
The concept of trimodule is a particular case of the concept of module over an algebra over an operad defined in \cite{[GK]}. For the more general context of $n$-ary algebras,  see \cite{[Gnedbaye]}.
\subsection{Trimodules and matched pairs of totally associative $3$-ary  algebras} 
\begin{definition}
A trimodule structure over totally associative $3$-ary algebra $(\A, \mu)$ on a vector space $V$ is defined by the following three linear multiplication mappings: 
\beqs
&&\mathcal{L}_{\mu}: \A\otimes\A\otimes V\rightarrow V, \cr
&& \mathcal{R}_{\mu}: V\otimes\A\otimes\A\rightarrow V, \cr
&& \mathcal{M}_{\mu}: \A\otimes V\otimes\A\rightarrow V
\eeqs 
 satisfying the following compatibility conditions
\beq \label{tr1}
\mathcal{L}_{\mu}(a, b)(\mathcal{L}_{\mu}(c, d)(v))= \mathcal{L}_{\mu}(\mu(a, b, c), d)(v)= \mathcal{L}_{\mu}(a, \mu(b, c, d))( v),
\eeq
\beq \label{tr2}
\mathcal{R}_{\mu}(c, d)(\mathcal{R}_{\mu}(a, b)(v))= \mathcal{R}_{\mu}(a, \mu(b, c, d))(v)= \mathcal{R}_{\mu}(\mu(a, b, c), d)(v),
\eeq
\beq \label{tr3}
\mathcal{M}_{\mu}(a, z)(\mathcal{M}_{\mu}(b, y)(\mathcal{M}_{\mu}(c, x)(v))= \mathcal{M}_{\mu}(\mu(a, b, c), \mu(x, y, z))(v),
\eeq 
\beq \label{tr4}
\mathcal{M}_{\mu}(a, d)(\mathcal{L}{\mu}(b, c)(v))= \mathcal{L}_{\mu}(a, b) (\mathcal{M}_{\mu}(c, d)(v))= \M_{\mu}(\mu(a, b, c), d)(v),
\eeq
\beq \label{tr5}
\mathcal{M}_{\mu}(a, d)(\mathcal{R}_{\mu}(b, c)(v))= \mathcal{R}_{\mu}(c, d)(\mathcal{M}_{\mu}(a, b)(v))= \mathcal{M}_{\mu}(a, \mu(b, c, d))(v),
\eeq
\beq \label{tr6}
\mathcal{R}_{\mu}(c, d)(\mathcal{L}_{\mu}(a, b) (v))= \mathcal{L}_{\mu}(a, b)(\mathcal{R}_{\mu}(c, d)(v))= \mathcal{M}_{\mu}(a, d)(\mathcal{M}_{\mu}(b, c)(v)),
\eeq
$ \forall a, b, c, d, x, y, z\in \A, v\in V$.
\end{definition} 

\begin{proposition}
($ \mathcal{L}_{\mu}, \mathcal{M}_{\mu}, \mathcal{R}_{\mu}, V$) is a trimodule of an associative totally $3$-ary algebra $(\mathcal{A}, \mu)$ if and only if the direct sum $(\mathcal{A} \oplus V, \tau) $ of the underlying vector spaces of $ \mathcal{A} $ and $V$ is turned into an associative   totally $3$-ary algebra $\tau$ given by
\beqs
\tau[(x + a), (y + b), (z + c)]= \mu(x, y, z) + \mathcal{L}_{\mu}(x, y)(c) + \mathcal{M}_{\mu}(x, z)(b) + \mathcal{R}_{\mu}(y, z)(a),
\eeqs
for all $ x, y, z \in \mathcal{A}, a, b, c \in V$. We denote it by $ \mathcal{A} \ltimes_{\mathcal{L}_{\mu}, \mathcal{M}_{\mu}, \mathcal{R}_{\mu}} V$.
\end{proposition}
\textbf{Proof}
Let $v_{1}, v_{2}, v_{3}, v_{4}, v_{5}\in V$ and $x_{1}, x_{2}, x_{3}, x_{4}, x_{5}\in \mathcal{A}$.
Set
\beqs
&&\tau[\tau[(x_{1} + v_{1}), (x_{2} + v_{2}), (x_{3} + v_{3})], (x_{4} + v_{4}), (x_{5} + v_{5})]\cr
&=&\tau[(x_{1} + v_{1}), \tau[(x_{2} + v_{2}), (x_{3} + v_{3}), (x_{4} + v_{4})], (x_{5} + v_{5})]\cr
&=&\tau[(x_{1} + v_{1}), (x_{2} + v_{2}), \tau[(x_{3} + v_{3}), (x_{4} + v_{4}), (x_{5} + v_{5})]]. 
\eeqs
After computation, we obtain the Eqs.(\ref{tr1})-(\ref{tr6}) while the Eq.(\ref{tr3}) is satisfied with the specification of the action of $\mathcal{M}_{\mu}$. Then ($ \mathcal{L}_{\mu}, \mathcal{M}_{\mu}, \mathcal{R}_{\mu}, V$) is a trimodule of the associative totally $3$-ary algebra $(\mathcal{A}, \mu)$ if and only if the direct sum $(\mathcal{A} \oplus V, \tau) $ is a totally $3$-ary algebra.

$ \hfill \square $
\begin{remark}
In the case where the Eq.(\ref{tr3}) is not satisfied, we refer to the name quasi trimodule structure instead of simply trimodule structure.
\end{remark}
\begin{theorem}\label{Theo. Match. Pair}
 Let $(\A, \mu_{\A})$  and $(\B, \mu_{\B})$ be two associative totally $3$-ary algebras. Suppose that
 there are linear maps $\mathcal{L}_{\mu_{\A}}: \A\otimes\A\otimes\B\rightarrow \B, \mathcal{R}_{\mu_{\A}}: \B\otimes\A\otimes\A\rightarrow \B, \mathcal{M}_{\mu_{\A}}: \A\otimes\B\otimes\A\rightarrow \B$ and $\mathcal{L}_{\mu_{\B}}: \B\otimes\B\otimes\A\rightarrow \A, \mathcal{R}_{\mu_{\B}}: \A\otimes\B\otimes\B\rightarrow \A, \mathcal{M}_{\mu_{\B}}: \B\otimes\A\otimes\B\rightarrow \A$ such that ($ \mathcal{L}_{\mu_{\A}}, \mathcal{M}_{\mu_{\A}}, \mathcal{R}_{\mu_{\A}}, \B$) is a quasi trimodule of the associative totally $3$-ary algebra $(\mathcal{A}, \mu_{A})$ and ($\mathcal{L}_{\mu_{\B}}, \mathcal{M}_{\mu_{\B}}, \mathcal{R}_{\mu_{\B}}, \A$) is a quasi trimodule of the associative totally $3$-ary algebra $(\mathcal{B}, \mu_{B})$ and they satisfy the following conditions:
 \beq \label{tr7}
\mu_{\A}(\mathcal{L}_{\mu_{\B}}(a, b)(x), y, z)= \mathcal{L}_{\mu_{\B}}[a, \mathcal{R}_{\mu_{\A}}(x, y)(b)](z)= \mathcal{L}_{\mu_{\B}}(a, b)(\mu_{\A}(x, y, z)),
\eeq
 \beq \label{tr8}
\mu_{\A}(\mathcal{M}_{\mu_{\B}}(a, b)(x), y, z)= \mathcal{L}_{\mu_{\B}}[a, \mathcal{M}_{\mu_{\A}}(x, y)(b)](z)= \mathcal{M}_{\mu_{\B}}[a, \mathcal{R}_{\mu_{\A}}(y, z)(b)](x),
\eeq
 \beq \label{tr9}
\mu_{\A}(\mathcal{R}_{\mu_{\B}}(a, b)(x), y, z)= \mu_{\A}[x, \mathcal{L}_{\mu_{\B}}(a, b)(y), z]= \mathcal{R}_{\mu_{\B}}[a, \mathcal{R}_{\mu_{\A}}(y, z)(b)](x),
\eeq
\beq \label{tr10}
\mathcal{L}_{\mu_{\B}}[\mathcal{L}_{\mu_{\A}}(x, y)(a), b](z)= \mu_{\A}[x, \mathcal{R}_{\mu_{\B}}(a, b)(y), z]= \mu_{\A}[x, y, \mathcal{L}_{\mu_{\B}}(a, b)(z)],
\eeq
\beq \label{tr11}
\mathcal{L}_{\mu_{\B}}[\mathcal{M}_{\mu_{\A}}(x, y)(a), b](z)= \mu_{\A}[x, \mathcal{M}_{\mu_{\B}}(a, b)(y), z]= \mathcal{R}_{\mu_{\B}}[a, \mathcal{M}_{\mu_{\A}}(y, z)(b)](x),
\eeq
\beq \label{tr12}
\mathcal{L}_{\mu_{\B}}[\mathcal{R}_{\mu_{\A}}(x, y)(a), b](z)= \mathcal{L}_{\mu_{\B}}[a, \mathcal{L}_{\mu_{\A}}(x, y)(b)](z)= \mathcal{M}_{\mu_{\B}}[a, \mathcal{M}_{\mu_{\A}}(y, z)(b)](x),
\eeq
\beq \label{tr13}
\mathcal{M}_{\mu_{\B}}[\mathcal{L}_{\mu_{\A}}(x, y)(a), b](z)= \mathcal{R}_{\mu_{\B}}[ \mathcal{M}_{\mu_{\A}}(y, z)(a), b](x)= \mu_{\A}[x, y, \mathcal{M}_{\mu_{\B}}(a, b)(z)],
\eeq
\beq \label{tr14}
\mathcal{M}_{\mu_{\B}}[\mathcal{M}_{\mu_{\A}}(x, y)(a), b](z)= \mathcal{R}_{\mu_{\B}}[ \mathcal{R}_{\mu_{\A}}(y, z)(a), b](x)= \mathcal{R}_{\mu_{\B}}[a, \mathcal{L}_{\mu_{\A}}(y, z)(b)](x),
\eeq
\beq \label{tr15}
\mathcal{M}_{\mu_{\B}}[\mathcal{R}_{\mu_{\A}}(x, y)(a), b](z)= \mathcal{M}_{\mu_{\B}}(a, b) (\mu_{\A}(x, y, z))= \mathcal{M}_{\mu_{\B}}[a, \mathcal{L}_{\mu_{\A}}(y, z)(b)](x),
\eeq
\beq \label{tr16}
\mathcal{R}_{\mu_{\B}}(a, b)(\mu_{\A}(x, y, z))= \mathcal{R}_{\mu_{\B}}[\mathcal{L}_{\mu_{\A}}(y, z)(a), b](x)= \mu_{\A}[x, y, \mathcal{R}_{\mu_{\B}}(a, b)(z)],
\eeq
\beq \label{tr17}
\mu_{\B}[\mathcal{L}_{\mu_{\A}}(x, y)(a), b, c]= \mathcal{L}_{\mu_{\A}}[x, \mathcal{R}_{\mu_{\B}}(a, b)(y)](c)= \mathcal{L}_{\mu_{\A}}(x, y)(\mu_{\B}(a, b, c)),
\eeq
\beq \label{tr18}
\mu_{\B}[\mathcal{M}_{\mu_{\A}}(x, y)(a), b, c]= \mathcal{L}_{\mu_{\A}}[x, \mathcal{M}_{\mu_{\B}}(a, b)(y)](c)= \mathcal{M}_{\mu_{\A}}[x, \mathcal{R}_{\mu_{\B}}(b, c)(y)](a),
\eeq
\beq \label{tr19}
\mu_{\B}[\mathcal{R}_{\mu_{\A}}(x, y)(a), b, c]= \mu_{\B}[a, \mathcal{L}_{\mu_{\A}}(x, y)(b), c]= \mathcal{R}_{\mu_{\A}}[x, \mathcal{R}_{\mu_{\B}}(b, c)(y)](a),
\eeq
\beq \label{tr20}
\mathcal{L}_{\mu_{\A}}[\mathcal{L}_{\mu_{\B}}(a, b)(x), y](c)= \mu_{\B}[a, \mathcal{R}_{\mu_{\A}}(x, y)(b), c]= \mu_{\B}[a, b, \mathcal{L}_{\mu_{\A}}(x, y)(c)],
\eeq
\beq \label{tr21}
\mathcal{L}_{\mu_{\A}}[\mathcal{M}_{\mu_{\B}}(a, b)(x), y](c)= \mu_{\B}[a, \mathcal{M}_{\mu_{\A}}(x, y)(b), c]= \mathcal{R}_{\mu_{\A}}[x, \mathcal{M}_{\mu_{\B}}(b, c)(y)](a),
\eeq
\beq \label{tr22}
\mathcal{L}_{\mu_{\A}}[\mathcal{R}_{\mu_{\B}}(a, b)(x), y](c)= \mathcal{L}_{\mu_{\A}}[x, \mathcal{L}_{\mu_{\B}}(a, b)(y)](c)= \mathcal{M}_{\mu_{\A}}[x, \mathcal{M}_{\mu_{\B}}(b, c)(y)](a),
\eeq
\beq \label{tr23}
\mathcal{M}_{\mu_{\A}}[\mathcal{L}_{\mu_{\B}}(a, b)(x), y](c)= \mathcal{R}_{\mu_{\A}}[ \mathcal{M}_{\mu_{\B}}(b, c)(x), y](a)= \mu_{\B}[a, b, \mathcal{M}_{\mu_{\A}}(x, y)(c)],
\eeq
\beq \label{tr24}
\mathcal{M}_{\mu_{\A}}[\mathcal{M}_{\mu_{\B}}(a, b)(x), y](c)= \mathcal{R}_{\mu_{\A}}[ \mathcal{R}_{\mu_{\B}}(b, c)(x), y](a)= \mathcal{R}_{\mu_{\A}}[x, \mathcal{L}_{\mu_{\B}}(b, c)(y)](a),
\eeq
\beq \label{tr25}
\mathcal{M}_{\mu_{\A}}[\mathcal{R}_{\mu_{\B}}(a, b)(x), y](c)= \mathcal{M}_{\mu_{\A}}(x, y)(\mu_{\B}(a, b, c))= \mathcal{M}_{\mu_{\A}}[x, \mathcal{L}_{\mu_{\B}}(b, c)(y)](a),
\eeq
\beq \label{tr26}
\mathcal{R}_{\mu_{\A}}(x, y)(\mu_{\B}(a, b, c))= \mathcal{R}_{\mu_{\A}}[\mathcal{L}_{\mu_{\B}}(b, c)(x), y](a)= \mu_{\B}[a, b, \mathcal{R}_{\mu_{\A}}(x, y)(c)],
\eeq
for any $ x, y, z\in \mathcal{A}, a, b, c\in \mathcal{B}$. Then, there is an associative totally $3$-ary  algebra structure on the direct sum $\mathcal{A} \oplus \mathcal{B}$ of the underlying vector spaces of $ \mathcal{A} $ and $ \mathcal{B} $ given by the product $\tau$ defined by 
\beqs
\tau[(x + a), (y + b), (z + c)]&=& [\mu_{\mathcal{A}}(x, y, z) + \mathcal{L}_{\mu_{\B}}(a, b)(z) + \mathcal{M}_{\mu_{\B}}(a, c)(y) + \mathcal{R}_{\mu_{\B}}(b, c)(x)] +\cr
&& [\mu_{\mathcal{B}}(a, b, c) + \mathcal{L}_{\mu_{\A}}(x, y)(c) + \mathcal{M}_{\mu_{\A}}(x, z)(b) + \mathcal{R}_{\mu_{\A}}(y, z)(a)]
\eeqs
for any $ x, y, z\in \mathcal{A}, a, b, c\in \mathcal{B}$. Let $ \mathcal{A} \bowtie^{\mathcal{L}_{\mu_{\mathcal{A}}}, \mathcal{M}_{\mu_{\mathcal{A}}}, 
\mathcal{R}_{\mu_{\mathcal{A}}}}_{\mathcal{L}_{\mu_{\mathcal{B}}}, \mathcal{M}_{\mu_{\mathcal{B}}}, \mathcal{R}_{\mu_{\mathcal{B}}}} \mathcal{B} $ denote this associative totally $3$-ary algebra.
\end{theorem}
\textbf{Proof}
Let $x_{1}, x_{2}, x_{3}, x_{4}, x_{5}\in \mathcal{A}$ and $y_{1}, y_{2}, y_{3}, y_{4}, y_{5}\in \mathcal{B}$. By definition, we have 
\beqs
\tau[(x + a), (y + b), (z + c)]&=& [\mu_{\mathcal{A}}(x, y, z) + \mathcal{L}_{\mu_{\B}}(a, b)(z) + \mathcal{M}_{\mu_{\B}}(a, c)(y) + \mathcal{R}_{\mu_{\B}}(b, c)(x)] +\cr
&& [\mu_{\mathcal{B}}(a, b, c) + \mathcal{L}_{\mu_{\A}}(x, y)(c) + \mathcal{M}_{\mu_{\A}}(x, z)(b) + \mathcal{R}_{\mu_{\A}}(y, z)(a)]
\eeqs 
for any $x, y, z\in \mathcal{A}, a, b, c\in \mathcal{B}$.
Setting the strong condition 
\beqs
&&\tau[\tau[(x_{1} + y_{1}), (x_{2} + y_{2}), (x_{3} + y_{3})], (x_{4} + y_{4}), (x_{5} + y_{5})]\cr
&=&\tau[(x_{1} + y_{1}), \tau[(x_{2} + y_{2}), (x_{3} + y_{3}), (x_{4} + y_{4})], (x_{5} + y_{5})]\cr
&=&\tau[(x_{1} + y_{1}), (x_{2} + y_{2}), \tau[(x_{3} + y_{3}), (x_{4} + y_{4}), (x_{5} + y_{5})]], 
\eeqs
we obtain by direct computation the Eqs.(\ref{tr7}) - (\ref{tr26}). Then, there is an associative totally $3$-ary  algebra structure on the direct sum $\mathcal{A} \oplus \mathcal{B}$ of the underlying vector spaces of $ \mathcal{A} $ and $ \mathcal{B} $ if and only if the Eqs.(\ref{tr7}) - (\ref{tr26}) are satisfied.
  
$ \hfill \square $
\begin{definition}
Let $ (\mathcal{A}, \mu_{\mathcal{A}}) $ and $  (\mathcal{B}, \mu_{\mathcal{B}}) $ 
be two totally associative $3$-ary algebras. Suppose that there are linear maps 
$\mathcal{L}_{\mu_{\mathcal{A}}}, \mathcal{M}_{\mu_{\mathcal{A}}}, \mathcal{R}_{\mu_{\mathcal{A}}}$ 
and $\mathcal{L}_{\mu_{\mathcal{B}}}, \mathcal{M}_{\mu_{\mathcal{B}}}, \mathcal{R}_{\mu_{\mathcal{B}}}$
 such that $(\mathcal{L}_{\mu_{\mathcal{A}}}, \mathcal{M}_{\mu_{\mathcal{A}}}, \mathcal{R}_{\mu_{\mathcal{A}}}, \mathcal{B})$ is a quasi trimodule of $ \mathcal{A} $ 
and $(\mathcal{L}_{\mu_{\mathcal{B}}}, \mathcal{M}_{\mu_{\mathcal{B}}}, \mathcal{R}_{\mu_{\mathcal{B}}}, \mathcal{A})$ is a quasi trimodule of $ \mathcal{B} $. 
If Eqs.\ref{tr7} - \ref{tr26} are satisfied, then $(\mathcal{A}, \mathcal{B},\mathcal{L}_{\mu_{\mathcal{A}}}, \mathcal{M}_{\mu_{\mathcal{A}}}, \mathcal{R}_{\mu_{\mathcal{A}}}, \mathcal{L}_{\mu_{\mathcal{B}}}, \mathcal{M}_{\mu_{\mathcal{B}}}, \mathcal{R}_{\mu_{\mathcal{B}}})$ 
 is called a matched pair of totally associative $3$-ary algebras.   
\end{definition}
\begin{definition}
Let $ (\mathcal{A}, \mu_{\mathcal{A}}) $ and $  (\mathcal{B}, \mu_{\mathcal{B}}) $ 
be two totally associative $3$-ary algebras. Suppose that there are linear maps 
$\mathcal{L}_{\mu_{\mathcal{A}}}, \mathcal{M}_{\mu_{\mathcal{A}}}, \mathcal{R}_{\mu_{\mathcal{A}}}$ 
and $\mathcal{L}_{\mu_{\mathcal{B}}}, \mathcal{M}_{\mu_{\mathcal{B}}}, \mathcal{R}_{\mu_{\mathcal{B}}}$
 such that $(\mathcal{L}_{\mu_{\mathcal{A}}}, \mathcal{M}_{\mu_{\mathcal{A}}}, \mathcal{R}_{\mu_{\mathcal{A}}}, \mathcal{B})$ is a trimodule of $ \mathcal{A} $ 
and $(\mathcal{L}_{\mu_{\mathcal{B}}}, \mathcal{M}_{\mu_{\mathcal{B}}}, \mathcal{R}_{\mu_{\mathcal{B}}}, \mathcal{A})$ is a trimodule of $ \mathcal{B} $. 
 
 Then $(\mathcal{A}, \mathcal{B},\mathcal{L}_{\mu_{\mathcal{A}}}, \mathcal{M}_{\mu_{\mathcal{A}}}, \mathcal{R}_{\mu_{\mathcal{A}}}, \mathcal{L}_{\mu_{\mathcal{B}}}, \mathcal{M}_{\mu_{\mathcal{B}}}, \mathcal{R}_{\mu_{\mathcal{B}}})$ 
 is a matched pair of totally associative $3$-ary algebras if Eqs.\ref{tr7} - \ref{tr26} are satisfied and  if the following conditions also  hold:
 \beq \label{tr27}
 \mathcal{M}_{\mu_{A}}(a, z)(\mathcal{M}_{\mu_{A}}(b, y)(\mathcal{M}_{\mu_{A}}(c, x)(v'))= \mathcal{M}_{\mu_{A}}(\mu_{A}(a, b, c), \mu_{A}(x, y, z))(v'),
 \eeq   
 \beq \label{tr28}
 \mathcal{M}_{\mu_{B}}(a', z')(\mathcal{M}_{\mu_{B}}(b', y')(\mathcal{M}_{\mu_{B}}(c', x')(v))= \mathcal{M}_{\mu_{B}}(\mu_{B}(a', b', c'), \mu_{B}(x', y', z'))(v),
 \eeq
 for any $ x, y, z, a, b, c, v\in \mathcal{A}, x', y', z', a', b', c', v'\in \mathcal{B}$.
\end{definition}
\begin{lemma}
Let $(\mathcal{L}_{\mu}, \mathcal{M}_{\mu}, \mathcal{R}_{\mu})$ be a trimodule of a totally associative $3$-ary algebra $\mathcal{A}$. Then, the linear maps   $\mathcal{L}^{\ast}_{\mu}, \mathcal{M}^{\ast}_{\mu}, \mathcal{R}^{\ast}_{\mu}: \mathcal{A}\otimes \mathcal{A\rightarrow} gl(V^{\ast})$ given by 
\beqs
&&\langle \mathcal{L}^{\ast}_{\mu}(x, y)u^{\ast}, v\rangle= \langle \mathcal{L}_{\mu}(x, y)v, u^{\ast}\rangle; \langle \mathcal{M}^{\ast}_{\mu}(x, y)u^{\ast}, v\rangle= \langle \mathcal{M}_{\mu}(x, y)v, u^{\ast}\rangle;\cr
&&\langle \mathcal{R}^{\ast}_{\mu}(x, y)u^{\ast}, v\rangle= \langle \mathcal{R}_{\mu}(x, y)v, u^{\ast}\rangle; \mbox{ for all } x, y \in \mathcal{A}, v \in V, u^{\ast}\in V^{\ast}
\eeqs
realize a trimodule of $\A$ denoted by $(\mathcal{R}^{\ast}_{\mu}, \mathcal{M}^{\ast}_{\mu}, \mathcal{L}^{\ast}_{\mu})$ satisfying 
\beq \label{dtr1}
\mathcal{R}^{\ast}_{\mu}(a, b)(\mathcal{R}^{\ast}_{\mu}(c, d)(u^{\ast}))= \mathcal{R}^{\ast}_{\mu}(\mu(a, b, c), d)(u^{\ast})= \mathcal{R}^{\ast}_{\mu}(a, \mu(b, c, d))(u^{\ast}),
\eeq
\beq \label{dtr2}
\mathcal{L}^{\ast}_{\mu}(c, d)(\mathcal{L}^{\ast}_{\mu}(a, b)(u^{\ast}))= \mathcal{L}^{\ast}_{\mu}(a, \mu(b, c, d))(u^{\ast})= \mathcal{L}^{\ast}_{\mu}(\mu(a, b, c), d)(u^{\ast}),
\eeq
\beq \label{dtr3}
\mathcal{M}^{\ast}_{\mu}(a, z)(\mathcal{M}^{\ast}_{\mu}(b, y)(\mathcal{M}^{\ast}_{\mu}(c, x)(u^{\ast}))= \mathcal{M}^{\ast}_{\mu}(\mu(c, b, a), \mu(z, y, x))(u^{\ast}),
\eeq 
\beq \label{dtr4}
\mathcal{M}^{\ast}_{\mu}(a, d)(\mathcal{R}^{\ast}_{\mu}{\mu}(b, c)(u^{\ast}))= \mathcal{R}^{\ast}_{\mu}(d, b) (\mathcal{M}^{\ast}_{\mu}(a, c)(u^{\ast}))= \M^{\ast}_{\mu}(a, \mu(d, b, c))(u^{\ast}),
\eeq
\beq \label{dtr5}
\mathcal{M}_{\mu}(a, d)(\mathcal{R}_{\mu}(b, c)(u^{\ast}))= \mathcal{L}^{\ast}_{\mu}(c, a)(\mathcal{M}^{\ast}_{\mu}(b, d)(u^{\ast}))= \mathcal{M}^{\ast}_{\mu}(\mu(b, c, a), d)(u^{\ast}),
\eeq
\beq \label{dtr6}
\mathcal{L}^{\ast}_{\mu}(c, d)(\mathcal{R}^{\ast}_{\mu}(a, b) (u^{\ast}))= \mathcal{R}^{\ast}_{\mu}(a, b)(\mathcal{L}^{\ast}_{\mu}(c, d)(u^{\ast}))= \mathcal{M}^{\ast}_{\mu}(d, a)(\mathcal{M}^{\ast}_{\mu}(c, b)(u^{\ast})),
\eeq
$ \forall a, b, c, d, x, y, z\in \A, u^{\ast}\in V^{\ast}$.
\end{lemma}
\textbf{Proof:}
By a direct computation, we obtain the result.

$ \hfill \square $

Let us now  give some notations useful in the sequel. Let $(\mathcal{A}, \mu)$ be an associative $3$-ary algebra.

 Considering the linear maps of the  left $L_{\mu},$ right $R_{\mu}$ and central $M_{\mu}$ multiplication operations defined as:
 \begin{eqnarray}
 L_{\mu}: \A\otimes \A & \longrightarrow & \mathfrak{gl}(\A)  \cr
  (x, y)  & \longmapsto & L_{\mu}(x, y):
  \begin{array}{ccc}
 \A &\longrightarrow & \A \cr 
  z & \longmapsto & \mu(x, y, z), 
   \end{array}
\end{eqnarray}
\begin{eqnarray}
    R_{\mu}: \A\otimes \A & \longrightarrow & \mathfrak{gl}(\A)  \cr
     (x, y)  & \longmapsto & R_{\mu}(x, y):
     \begin{array}{ccc}
    \A &\longrightarrow & \A \cr 
     z & \longmapsto & \mu(z, x, y),
      \end{array}
\end{eqnarray}
      \begin{eqnarray}
    M_{\mu}: \A\otimes \A & \longrightarrow & \mathfrak{gl}(\A)  \cr
     (x, y)  & \longmapsto & M_{\mu}(x, y):
     \begin{array}{ccc}
    \A &\longrightarrow & \A \cr 
     z & \longmapsto & \mu(x, z, y).
      \end{array}
 \end{eqnarray}
 The dual maps $L^{*}_{\mu}, R^{*}_{\mu}, M^{*}_{\mu}$  of the linear maps $L_{\mu}, R_{\mu}, M_{\mu}$ are defined, respectively, as:\\ 
$\displaystyle L^{*}_{\mu}, R^{*}_{\mu}, M^{*}_{\mu}: \A\otimes \A \rightarrow \mathfrak{gl}(\A^{*})$
 such that:
\beq\label{dual1}
 L^*_{\mu}: \A\otimes \A & \longrightarrow & \mathfrak{gl}(\A^*)  \cr
  (x, y)  & \longmapsto & L^*_\mu(x, y):
      \begin{array}{llll}
 \A^* &\longrightarrow & \A^* \\ 
  u^* & \longmapsto & L^*_\mu(x, y)(u^*): 
      \begin{array}{llll}
\A  \longrightarrow   \mathcal{K} \cr
v  \longmapsto  &\left< L^{*}_\mu(x, y)(u^{*}), v\right>\cr
 &:= \left<  u^{*}, L_\mu(x, y) v\right>, 
\end{array} 
     \end{array}
 \eeq
\beq\label{dual2}
 R^*_{\mu}: \A\otimes \A & \longrightarrow & \mathfrak{gl}(\A^*)  \cr
  (x, y)  & \longmapsto & R^*_\mu(x, y):
      \begin{array}{llll}
 \A^* &\longrightarrow & \A^* \\ 
  u^* & \longmapsto & R^*_\mu(x, y)(u^*): 
      \begin{array}{llll}
\A  \longrightarrow  \mathcal{K} \cr
v  \longmapsto & \left< R^{*}_\mu(x, y)(u^{*}), v\right>\cr
 &:= \left<  u^{*}, R_\mu(x, y) v\right>, 
\end{array} 
     \end{array}
     \eeq
 
  \beq\label{dual3}
 M^*_{\mu}: \A\otimes \A & \longrightarrow & \mathfrak{gl}(\A^*)  \cr
  (x, y)  & \longmapsto & M^*_\mu(x, y):
      \begin{array}{llll}
 \A^* &\longrightarrow & \A^* \\ 
  u^* & \longmapsto & M^*_\mu(x, y)(u^*): 
      \begin{array}{llll}
\A  \longrightarrow  \mathcal{K} \cr
v\longmapsto &\left< M^{*}_\mu(x, y)(u^{*}), v\right>\cr
 &:= \left<  u^{*}, M_\mu(x, y) v\right>, 
\end{array} 
     \end{array}
 \eeq
for all $x, y, z, v \in \A, u^{*} \in \A^{*},$ where $\A^{*}$ is the dual space of $\A.$ 
\begin{proposition}
Let $(\mathcal{A}, \mu)$ be a totally associative $3$-ary algebra. Then $(L_{\mu}, 0, 0,
 \mathcal{A})$,\\ $(0, 0, R_{\mu}, \mathcal{A})$, $(L_{\mu}, M_{\mu}, R_{\mu}, 
 \mathcal{A})$, $(R^{\ast}_{\mu}, 0, 0, \mathcal{A}^{\ast})$, $(0, 0, 
 L^{\ast}_{\mu}, \mathcal{A}^{\ast})$ and $(R^{\ast}_{\mu}, M^{\ast}_{\mu}, L^{\ast}_{\mu}, 
 \mathcal{A}^{\ast})$ are  trimodules of $\mathcal{A}$.
\end{proposition}
\textbf{Proof:}
It  is  straightforward.

$ \hfill \square $

 For a linear map $\phi: V_{1}\rightarrow V_{2}$, we denote the dual (linear) map by $\phi^{\ast}: V^{\ast}_{2}\rightarrow V^{\ast}_{1}$ given by
 \beqs
 \langle v, \phi^{\ast}(u^{\ast})\rangle= \langle \phi(v), u^{\ast}\rangle \mbox{ for all } v\in V_{1}, u^{\ast}\in V^{\ast}_{2}.
 \eeqs


\begin{theorem}
 Let $(\mathcal{A}, \mu)$  be a totally associative $3$-ary algebra. Suppose that there is a totally associative $3$-ary algebra structure $\nu$ on its dual space $\mathcal{A}^{\ast}$. Then, $(\mathcal{A}, \mathcal{A}^{\ast} , R^{\ast}_{\mu} , M^{\ast}_{\mu} , L^{\ast}_{\mu}, R^{\ast}_{\nu} , M^{\ast}_{\nu} , L^{\ast}_{\nu})$ is a matched pair of totally associative $3$-ary algebras if and only if there is an associative totally $3$-ary  algebra structure on the direct sum $\mathcal{A} \oplus \mathcal{A}^{\ast}$ of the underlying vector spaces of $ \mathcal{A} $ and $ \mathcal{A}^{\ast} $ given by the product $\tau$ defined by
 \beqs
\tau[(x + a^{\ast}), (y + b^{\ast}), (z + c^{\ast})]&=& [\mu(x, y, z) + L^{\ast}_{\nu}(a^{\ast}, b^{\ast})(z) + M^{\ast}_{\nu}(a^{\ast}, c^{\ast})(y) + R^{\ast}_{\nu}(b^{\ast}, c^{\ast})(x)] +\cr
&& [\nu(a^{\ast}, b^{\ast}, c^{\ast}) + L^{\ast}_{\mu}(x, y)(c^{\ast}) + M^{\ast}_{\mu}(x, z)(b^{\ast}) + R^{\ast}_{\mu}(y, z)(a^{\ast})]
\eeqs
for any $x, y, z\in \mathcal{A}, a^{\ast}, b^{\ast}, c^{\ast}\in \mathcal{A}^{\ast}$.
\end{theorem}
\textbf{Proof}
We hold a same reasoning as in the proof of Theorem(\ref{Theo. Match. Pair}) with $\mathcal{B}= \mathcal{A}^{\ast}$ and
\beqs
\tau[(x + a^{\ast}), (y + b^{\ast}), (z + c^{\ast})]&=& [\mu(x, y, z) + L^{\ast}_{\nu}(a^{\ast}, b^{\ast})(z) + M^{\ast}_{\nu}(a^{\ast}, c^{\ast})(y) + R^{\ast}_{\nu}(b^{\ast}, c^{\ast})(x)] +\cr
&& [\nu(a^{\ast}, b^{\ast}, c^{\ast}) + L^{\ast}_{\mu}(x, y)(c^{\ast}) + M^{\ast}_{\mu}(x, z)(b^{\ast}) + R^{\ast}_{\mu}(y, z)(a^{\ast})]
\eeqs
for any $x, y, z\in \mathcal{A}, a^{\ast}, b^{\ast}, c^{\ast}\in \mathcal{A}^{\ast}$.
$ \hfill \square $

\subsection{Trimodules and matched pairs of partially associative $3$-ary  algebras}

In this sequel, we give a definition of trimodule for the partially associative $3$-ary algebras. Then, the matched pairs of partially associative $3$-ary algebras are established. 
\begin{definition}\label{Def4}
A trimodule structure over partially associative $3$-ary algebra $(\A, \mu)$ on a vector space $V$ is defined by the following three linear multiplication mappings: 
\beqs
&&\mathcal{L}_{\mu}: \A\otimes\A\otimes V\rightarrow V, \cr
&& \mathcal{R}_{\mu}: V\otimes\A\otimes\A\rightarrow V, \cr
&& \mathcal{M}_{\mu}: \A\otimes V\otimes\A\rightarrow V
\eeqs 
 satisfying the following compatibility conditions
\beq \label{p1}
\mathcal{L}_{\mu}(a, b)(\mathcal{L}_{\mu}(c, d)(v)) + \mathcal{L}_{\mu}(\mu(a, b, c), d)(v) + \mathcal{L}_{\mu}(a, \mu(b, c, d))( v)=0,
\eeq
\beq \label{p2}
\mathcal{R}_{\mu}(c, d)(\mathcal{R}_{\mu}(a, b)(v)) + \mathcal{R}_{\mu}(a, \mu(b, c, d))(v) + \mathcal{R}_{\mu}(\mu(a, b, c), d)(v)=0,
\eeq
\beq \label{p3}
\mathcal{M}_{\mu}(a, z)(\mathcal{M}_{\mu}(b, y)(\mathcal{M}_{\mu}(c, x)(v)) = \mathcal{M}_{\mu}(\mu(a, b, c), \mu(x, y, z))(v),
\eeq 
\beq \label{p4}
\mathcal{M}_{\mu}(a, d)(\mathcal{L}{\mu}(b, c)(v)) + \mathcal{L}_{\mu}(a, b) (\mathcal{M}_{\mu}(c, d)(v)) + \M_{\mu}(\mu(a, b, c), d)(v)=0,
\eeq
\beq \label{p5}
\mathcal{M}_{\mu}(a, d)(\mathcal{R}_{\mu}(b, c)(v)) + \mathcal{R}_{\mu}(c, d)(\mathcal{M}_{\mu}(a, b)(v)) + \mathcal{M}_{\mu}(a, \mu(b, c, d))(v)=0,
\eeq
\beq \label{p6}
\mathcal{R}_{\mu}(c, d)(\mathcal{L}_{\mu}(a, b) (v)) + \mathcal{L}_{\mu}(a, b)(\mathcal{R}_{\mu}(c, d)(v)) + \mathcal{M}_{\mu}(a, d)(\mathcal{M}_{\mu}(b, c)(v))=0,
\eeq
$ \forall a, b, c, d, x, y, z\in \A, v\in V$.
\end{definition} 
\begin{proposition}
($ \mathcal{L}_{\mu}, \mathcal{M}_{\mu}, \mathcal{R}_{\mu}, V$) is a trimodule of an associative partially $3$-ary algebra $(\mathcal{A}, \mu)$ if and only if the direct sum $(\mathcal{A} \oplus V, \tau) $ of the underlying vector spaces of $ \mathcal{A} $ and $V$ is turned into an associative   partially $3$-ary algebra $\tau$ given by
\beqs
\tau[(x + a), (y + b), (z + c)]= \mu(x, y, z) + \mathcal{L}_{\mu}(x, y)(c) + \mathcal{M}_{\mu}(x, z)(b) + \mathcal{R}_{\mu}(y, z)(a),
\eeqs
for all $ x, y, z \in \mathcal{A}, a, b, c \in V$. We denote it by $ \mathcal{A} \ltimes_{\mathcal{L}_{\mu}, \mathcal{M}_{\mu}, \mathcal{R}_{\mu}} V$.
\end{proposition}
\textbf{Proof}
Let $v_{1}, v_{2}, v_{3}, v_{4}, v_{5}\in V$ and $x_{1}, x_{2}, x_{3}, x_{4}, x_{5}\in \mathcal{A}$.
Set
\beqs
&&\tau[\tau[(x_{1} + v_{1}), (x_{2} + v_{2}), (x_{3} + v_{3})], (x_{4} + v_{4}), (x_{5} + v_{5})]+\cr
&&\tau[(x_{1} + v_{1}), \tau[(x_{2} + v_{2}), (x_{3} + v_{3}), (x_{4} + v_{4})], (x_{5} + v_{5})]+\cr
&&\tau[(x_{1} + v_{1}), (x_{2} + v_{2}), \tau[(x_{3} + v_{3}), (x_{4} + v_{4}), (x_{5} + v_{5})]]=0. 
\eeqs
After computation, we obtain the Eqs.(\ref{p1})-(\ref{p6}) while the Eq.(\ref{p3}) is satisfied with the specification of the action of $\mathcal{M}_{\mu}$. Then ($ \mathcal{L}_{\mu}, \mathcal{M}_{\mu}, \mathcal{R}_{\mu}, V$) is a trimodule of the associative partially $3$-ary algebra $(\mathcal{A}, \mu)$ if and only if the direct sum $(\mathcal{A} \oplus V, \tau) $ is a partially $3$-ary algebra.

$ \hfill \square $
\begin{remark}
In the case where the Eq.(\ref{p3}) is not satisfied, we refer to the name of quasi trimodule structure instead of simply trimodule structure.
\end{remark}
\begin{theorem}\label{Theo. Match. Partially}
 Let $(\A, \mu_{\A})$  and $(\B, \mu_{\B})$ be two associative partially $3$-ary algebras. Suppose that
 there are linear maps $\mathcal{L}_{\mu_{\A}}: \A\otimes\A\otimes\B\rightarrow \B, \mathcal{R}_{\mu_{\A}}: \B\otimes\A\otimes\A\rightarrow \B, \mathcal{M}_{\mu_{\A}}: \A\otimes\B\otimes\A\rightarrow \B$ and $\mathcal{L}_{\mu_{\B}}: \B\otimes\B\otimes\A\rightarrow \A, \mathcal{R}_{\mu_{\B}}: \A\otimes\B\otimes\B\rightarrow \A, \mathcal{M}_{\mu_{\B}}: \B\otimes\A\otimes\B\rightarrow \A$ such that ($ \mathcal{L}_{\mu_{\A}}, \mathcal{M}_{\mu_{\A}}, \mathcal{R}_{\mu_{\A}}, \B$) is a quasi trimodule of the associative partially $3$-ary algebra $(\mathcal{A}, \mu_{A}),$ and ($ \mathcal{L}_{\mu_{\B}}, \mathcal{M}_{\mu_{\B}}, \mathcal{R}_{\mu_{\B}}, \A$) is a quasi trimodule of the associative partially $3$-ary algebra $(\mathcal{B}, \mu_{B}),$  satisfying the following conditions:
 \beq \label{p7}
\mu_{\A}(\mathcal{L}_{\mu_{\B}}(a, b)(x), y, z) + \mathcal{L}_{\mu_{\B}}[a, \mathcal{R}_{\mu_{\A}}(y, z)(b)](z) + \mathcal{L}_{\mu_{\B}}(a, b)(\mu_{\A}(x, y, z))=0,
\eeq
 \beq \label{p8}
\mu_{\A}(\mathcal{M}_{\mu_{\B}}(a, b)(x), y, z) + \mathcal{L}_{\mu_{\B}}[a, \mathcal{M}_{\mu_{\A}}(x, y)(b)](z) + \mathcal{M}_{\mu_{\B}}[a, \mathcal{R}_{\mu_{\A}}(y, z)(b)](x)=0,
\eeq
 \beq \label{p9}
\mu_{\A}(\mathcal{R}_{\mu_{\B}}(a, b)(x), y, z) + \mu_{\A}[x, \mathcal{L}_{\mu_{\B}}(a, b)(y), z] + \mathcal{R}_{\mu_{\B}}[a, \mathcal{R}_{\mu_{\A}}(y, z)(b)](x)=0,
\eeq
\beq \label{p10}
\mathcal{L}_{\mu_{\B}}[\mathcal{L}_{\mu_{\A}}(x, y)(a), b](z) + \mu_{\A}[x, \mathcal{R}_{\mu_{\B}}(a, b)(y), z] + \mu_{\A}[x, y, \mathcal{L}_{\mu_{\B}}(a, b)(z)]=0,
\eeq
\beq \label{p11}
\mathcal{L}_{\mu_{\B}}[\mathcal{M}_{\mu_{\A}}(x, y)(a), b](z) + \mu_{\A}[x, \mathcal{M}_{\mu_{\B}}(a, b)(y), z] + \mathcal{R}_{\mu_{\B}}[a, \mathcal{M}_{\mu_{\A}}(y, z)(b)](x)=0,
\eeq
\beq \label{p12}
\mathcal{L}_{\mu_{\B}}[\mathcal{R}_{\mu_{\A}}(x, y)(a), b](z) + \mathcal{L}_{\mu_{\B}}[a, \mathcal{L}_{\mu_{\A}}(x, y)(b)](z) + \mathcal{M}_{\mu_{\B}}[a, \mathcal{M}_{\mu_{\A}}(y, z)(b)](x)=0,
\eeq
\beq \label{p13}
\mathcal{M}_{\mu_{\B}}[\mathcal{L}_{\mu_{\A}}(x, y)(a), b](z) + \mathcal{R}_{\mu_{\B}}[ \mathcal{M}_{\mu_{\A}}(y, z)(a), b](x) + \mu_{\A}[x, y, \mathcal{M}_{\mu_{\B}}(a, b)(z)]=0,
\eeq
\beq \label{p14}
\mathcal{M}_{\mu_{\B}}[\mathcal{M}_{\mu_{\A}}(x, y)(a), b](z) + \mathcal{R}_{\mu_{\B}}[ \mathcal{R}_{\mu_{\A}}(y, z)(a), b](x) + \mathcal{R}_{\mu_{\B}}[a, \mathcal{L}_{\mu_{\A}}(y, z)(b)](x)=0,
\eeq
\beq \label{p15}
\mathcal{M}_{\mu_{\B}}[\mathcal{R}_{\mu_{\A}}(x, y)(a), b](z) + \mathcal{M}_{\mu_{\B}}(a, b) (\mu_{\A}(x, y, z) + \mathcal{M}_{\mu_{\B}}[a, \mathcal{L}_{\mu_{\A}}(y, z)(b)](x)=0,
\eeq
\beq \label{p16}
\mathcal{R}_{\mu_{\B}}(a, b)(\mu_{\A}(x, y, z)) + \mathcal{R}_{\mu_{\B}}[\mathcal{L}_{\mu_{\A}}(y, z)(a), b](x) + \mu_{\A}[x, y, \mathcal{R}_{\mu_{\B}}(a, b)(z)]=0,
\eeq
\beq \label{p17}
\mu_{\B}[\mathcal{L}_{\mu_{\A}}(x, y)(a), b, c] + \mathcal{L}_{\mu_{\A}}[x, \mathcal{R}_{\mu_{\B}}(a, b)(y)](c) + \mathcal{L}_{\mu_{\A}}(x, y)(\mu_{\B}(a, b, c))=0,
\eeq
\beq \label{p18}
\mu_{\B}[\mathcal{M}_{\mu_{\A}}(x, y)(a), b, c] + \mathcal{L}_{\mu_{\A}}[x, \mathcal{M}_{\mu_{\B}}(a, b)(y)](c) + \mathcal{M}_{\mu_{\A}}[x, \mathcal{R}_{\mu_{\B}}(b, c)(y)](a)=0,
\eeq
\beq \label{p19}
\mu_{\B}[\mathcal{R}_{\mu_{\A}}(x, y)(a), b, c] + \mu_{\B}[a, \mathcal{L}_{\mu_{\A}}(x, y)(b), c] + \mathcal{R}_{\mu_{\A}}[x, \mathcal{R}_{\mu_{\B}}(b, c)(y)](a)=0,
\eeq
\beq \label{p20}
\mathcal{L}_{\mu_{\A}}[\mathcal{L}_{\mu_{\B}}(a, b)(x), y](c) + \mu_{\B}[a, \mathcal{R}_{\mu_{\A}}(x, y)(b), c] + \mu_{\B}[a, b, \mathcal{L}_{\mu_{\A}}(x, y)(c)]=0,
\eeq
\beq \label{p21}
\mathcal{L}_{\mu_{\A}}[\mathcal{M}_{\mu_{\B}}(a, b)(x), y](c) + \mu_{\B}[a, \mathcal{M}_{\mu_{\A}}(x, y)(b), c] + \mathcal{R}_{\mu_{\A}}[x, \mathcal{M}_{\mu_{\B}}(b, c)(y)](a)=0,
\eeq
\beq \label{p22}
\mathcal{L}_{\mu_{\A}}[\mathcal{R}_{\mu_{\B}}(a, b)(x), y](c) + \mathcal{L}_{\mu_{\A}}[x, \mathcal{L}_{\mu_{\B}}(a, b)(y)](c) + \mathcal{M}_{\mu_{\A}}[x, \mathcal{M}_{\mu_{\B}}(b, c)(y)](a)=0,
\eeq
\beq \label{p23}
\mathcal{M}_{\mu_{\A}}[\mathcal{L}_{\mu_{\B}}(a, b)(x), y](c) + \mathcal{R}_{\mu_{\A}}[ \mathcal{M}_{\mu_{\B}}(b, c)(x), y](a) + \mu_{\B}[a, b, \mathcal{M}_{\mu_{\A}}(x, y)(c)]=0,
\eeq
\beq \label{p24}
\mathcal{M}_{\mu_{\A}}[\mathcal{M}_{\mu_{\B}}(a, b)(x), y](c) +\mathcal{R}_{\mu_{\A}}[ \mathcal{R}_{\mu_{\B}}(b, c)(x), y](a) + \mathcal{R}_{\mu_{\A}}[x, \mathcal{L}_{\mu_{\B}}(b, c)(y)](a)=0,
\eeq
\beq \label{p25}
\mathcal{M}_{\mu_{\A}}[\mathcal{R}_{\mu_{\B}}(a, b)(x), y](c) + \mathcal{M}_{\mu_{\A}}(x, y)(\mu_{\B}(a, b, c)) + \mathcal{M}_{\mu_{\A}}[x, \mathcal{L}_{\mu_{\B}}(b, c)(y)](a)=0,
\eeq
\beq \label{p26}
\mathcal{R}_{\mu_{\A}}(x, y)(\mu_{\B}(a, b, c)) + \mathcal{R}_{\mu_{\A}}[\mathcal{L}_{\mu_{\B}}(b, c)(x), y](a) + \mu_{\B}[a, b, \mathcal{R}_{\mu_{\A}}(x, y)(c)]=0,
\eeq
for any $ x, y, z\in \mathcal{A}, a, b, c\in \mathcal{B}$. Then, there is an associative partially $3$-ary  algebra structure on the direct sum $\mathcal{A} \oplus \mathcal{B}$ of the underlying vector spaces of $ \mathcal{A} $ and $ \mathcal{B} $ given by the product $\tau$ defined by 
\beqs
\tau[(x + a), (y + b), (z + c)]&=& [\mu_{\mathcal{A}}(x, y, z) + \mathcal{L}_{\mu_{\B}}(a, b)(z) + \mathcal{M}_{\mu_{\B}}(a, c)(y) + \mathcal{R}_{\mu_{\B}}(b, c)(x)] +\cr
&& [\mu_{\mathcal{B}}(a, b, c) + \mathcal{L}_{\mu_{\A}}(x, y)(c) + \mathcal{M}_{\mu_{\A}}(x, z)(b) + \mathcal{R}_{\mu_{\A}}(y, z)(a)]
\eeqs
for any $ x, y, z\in \mathcal{A}, a, b, c\in \mathcal{B}$. Let $ \mathcal{A} \bowtie^{\mathcal{L}_{\mu_{\mathcal{A}}}, \mathcal{M}_{\mu_{\mathcal{A}}}, 
\mathcal{R}_{\mu_{\mathcal{A}}}}_{\mathcal{L}_{\mu_{\mathcal{B}}}, \mathcal{M}_{\mu_{\mathcal{B}}}, \mathcal{R}_{\mu_{\mathcal{B}}}} \mathcal{B} $ denote this associative partially $3$-ary algebra.
\end{theorem}
\textbf{Proof}
Let $x_{1}, x_{2}, x_{3}, x_{4}, x_{5}\in \mathcal{A}$ and $y_{1}, y_{2}, y_{3}, y_{4}, y_{5}\in \mathcal{B}$. By definition, we have 
\beqs
\tau[(x + a), (y + b), (z + c)]&=& [\mu_{\mathcal{A}}(x, y, z) + \mathcal{L}_{\mu_{\B}}(x, y)(c) + \mathcal{M}_{\mu_{\B}}(x, z)(b) + \mathcal{R}_{\mu_{\B}}(y, z)(a)] +\cr
&& [\mu_{\mathcal{B}}(a, b, c) + \mathcal{L}_{\mu_{\A}}(a, b)(z) + \mathcal{M}_{\mu_{\A}}(a, c)(y) + \mathcal{R}_{\mu_{\A}}(b, c)(x)]
\eeqs 
for any $x, y, z\in \mathcal{A}, a, b, c\in \mathcal{B}$.
Setting the strong condition 
\beqs
&&\tau[\tau[(x_{1} + y_{1}), (x_{2} + y_{2}), (x_{3} + y_{3})], (x_{4} + y_{4}), (x_{5} + y_{5})]\cr
&=&\tau[(x_{1} + y_{1}), \tau[(x_{2} + y_{2}), (x_{3} + y_{3}), (x_{4} + y_{4})], (x_{5} + y_{5})]\cr
&=&\tau[(x_{1} + y_{1}), (x_{2} + y_{2}), \tau[(x_{3} + y_{3}), (x_{4} + y_{4}), (x_{5} + y_{5})]], 
\eeqs
we obtain by direct computation the Eqs.(\ref{p7}) - (\ref{p26}). Then, there is an associative partially $3$-ary  algebra structure on the direct sum $\mathcal{A} \oplus \mathcal{B}$ of the underlying vector spaces of $ \mathcal{A} $ and $ \mathcal{B} $ if and only if the Eqs.(\ref{p7}) - (\ref{p26}) are satisfied.
  
$ \hfill \square $
\begin{definition}
Let $ (\mathcal{A}, \mu_{\mathcal{A}}) $ and $  (\mathcal{B}, \mu_{\mathcal{B}}) $ 
be two partially associative $3$-ary algebras. Suppose that there are linear maps 
$\mathcal{L}_{\mu_{\mathcal{A}}}, \mathcal{M}_{\mu_{\mathcal{A}}}, \mathcal{R}_{\mu_{\mathcal{A}}}$ 
and $\mathcal{L}_{\mu_{\mathcal{B}}}, \mathcal{M}_{\mu_{\mathcal{B}}}, \mathcal{R}_{\mu_{\mathcal{B}}}$
 such that $(\mathcal{L}_{\mu_{\mathcal{A}}}, \mathcal{M}_{\mu_{\mathcal{A}}}, \mathcal{R}_{\mu_{\mathcal{A}}}, \mathcal{B})$ is a quasi trimodule of $ \mathcal{A}, $ 
and $(\mathcal{L}_{\mu_{\mathcal{B}}}, \mathcal{M}_{\mu_{\mathcal{B}}}, \mathcal{R}_{\mu_{\mathcal{B}}}, \mathcal{A})$ is a quasi trimodule of $ \mathcal{B} $. 
If Eqs.\ref{p7} - \ref{p26} are satisfied, then $(\mathcal{A}, \mathcal{B},\mathcal{L}_{\mu_{\mathcal{A}}}, \mathcal{M}_{\mu_{\mathcal{A}}}, \mathcal{R}_{\mu_{\mathcal{A}}}, \mathcal{L}_{\mu_{\mathcal{B}}}, \mathcal{M}_{\mu_{\mathcal{B}}}, \mathcal{R}_{\mu_{\mathcal{B}}})$ 
 is called a matched pair of partially associative $3$-ary algebras.   
\end{definition}
\begin{definition}
Let $ (\mathcal{A}, \mu_{\mathcal{A}}) $ and $  (\mathcal{B}, \mu_{\mathcal{B}}) $ 
be two partially associative $3$-ary algebras. Suppose that there are linear maps 
$\mathcal{L}_{\mu_{\mathcal{A}}}, \mathcal{M}_{\mu_{\mathcal{A}}}, \mathcal{R}_{\mu_{\mathcal{A}}}$ 
and $\mathcal{L}_{\mu_{\mathcal{B}}}, \mathcal{M}_{\mu_{\mathcal{B}}}, \mathcal{R}_{\mu_{\mathcal{B}}}$
 such that $(\mathcal{L}_{\mu_{\mathcal{A}}}, \mathcal{M}_{\mu_{\mathcal{A}}}, \mathcal{R}_{\mu_{\mathcal{A}}}, \mathcal{B})$ is a trimodule of $ \mathcal{A}, $ 
and $(\mathcal{L}_{\mu_{\mathcal{B}}}, \mathcal{M}_{\mu_{\mathcal{B}}}, \mathcal{R}_{\mu_{\mathcal{B}}}, \mathcal{A})$ is a trimodule of $ \mathcal{B} $. 
 
 Then $(\mathcal{A}, \mathcal{B},\mathcal{L}_{\mu_{\mathcal{A}}}, \mathcal{M}_{\mu_{\mathcal{A}}}, \mathcal{R}_{\mu_{\mathcal{A}}}, \mathcal{L}_{\mu_{\mathcal{B}}}, \mathcal{M}_{\mu_{\mathcal{B}}}, \mathcal{R}_{\mu_{\mathcal{B}}})$ 
 is a matched pair of partially associative $3$-ary algebras if Eqs.\ref{p7} - \ref{p26} and the following two conditions are satisfied:
 \beq \label{p27}
 \mathcal{M}_{\mu_{A}}(a, z)(\mathcal{M}_{\mu_{A}}(b, y)(\mathcal{M}_{\mu_{A}}(c, x)(v'))= \mathcal{M}_{\mu_{A}}(\mu_{A}(a, b, c), \mu_{A}(x, y, z))(v'),
 \eeq   
 \beq \label{p28}
 \mathcal{M}_{\mu_{B}}(a', z')(\mathcal{M}_{\mu_{B}}(b', y')(\mathcal{M}_{\mu_{B}}(c', x')(v))= \mathcal{M}_{\mu_{B}}(\mu_{B}(a', b', c'), \mu_{B}(x', y', z'))(v),
 \eeq
 for any $ x, y, z, a, b, c, v\in \mathcal{A}, x', y', z', a', b', c', v'\in \mathcal{B}$.
\end{definition}
\begin{lemma}
Let $(\mathcal{L}_{\mu}, \mathcal{M}_{\mu}, \mathcal{R}_{\mu})$ be a trimodule of a partially associative $3$-ary algebra $\mathcal{A}$. Then, the linear maps   $\mathcal{L}^{\ast}_{\mu}, \mathcal{M}^{\ast}_{\mu}, \mathcal{R}^{\ast}_{\mu}: \mathcal{A}\otimes \mathcal{A\rightarrow} gl(V^{\ast})$ given by 
\beqs
&&\langle \mathcal{L}^{\ast}_{\mu}(x, y)u^{\ast}, v\rangle= \langle \mathcal{L}_{\mu}(x, y)v, u^{\ast}\rangle; \langle \mathcal{M}^{\ast}_{\mu}(x, y)u^{\ast}, v\rangle= \langle \mathcal{M}_{\mu}(x, y)v, u^{\ast}\rangle;\cr
&&\langle \mathcal{R}^{\ast}_{\mu}(x, y)u^{\ast}, v\rangle= \langle \mathcal{R}_{\mu}(x, y)v, u^{\ast}\rangle; \mbox{ for all } x, y \in \mathcal{A}, v \in V, u^{\ast}\in V^{\ast}.
\eeqs
realize a trimodule of $\A$ denoted by $(\mathcal{R}^{\ast}_{\mu}, \mathcal{M}^{\ast}_{\mu}, \mathcal{L}^{\ast}_{\mu})$ satisfying 
\beq \label{dpr1}
\mathcal{R}^{\ast}_{\mu}(a, b)(\mathcal{R}^{\ast}_{\mu}(c, d)(u^{\ast})) + \mathcal{R}^{\ast}_{\mu}(\mu(a, b, c), d)(u^{\ast}) + \mathcal{R}^{\ast}_{\mu}(a, \mu(b, c, d))(u^{\ast})= 0,
\eeq
\beq \label{dpr2}
\mathcal{L}^{\ast}_{\mu}(c, d)(\mathcal{L}^{\ast}_{\mu}(a, b)(u^{\ast})) + \mathcal{L}^{\ast}_{\mu}(a, \mu(b, c, d))(u^{\ast}) + \mathcal{L}^{\ast}_{\mu}(\mu(a, b, c), d)(u^{\ast})= 0,
\eeq
\beq \label{dpr3}
\mathcal{M}^{\ast}_{\mu}(a, z)(\mathcal{M}^{\ast}_{\mu}(b, y)(\mathcal{M}^{\ast}_{\mu}(c, x)(u^{\ast}))=\mathcal{M}^{\ast}_{\mu}(\mu(c, b, a), \mu(z, y, x))(u^{\ast}),
\eeq 
\beq \label{dtr4}
\mathcal{M}^{\ast}_{\mu}(a, d)(\mathcal{R}^{\ast}_{\mu}{\mu}(b, c)(u^{\ast})) + \mathcal{R}^{\ast}_{\mu}(d, b) (\mathcal{M}^{\ast}_{\mu}(a, c)(u^{\ast})) + \M^{\ast}_{\mu}(a, \mu(d, b, c))(u^{\ast})= 0,
\eeq
\beq \label{dtr5}
\mathcal{M}_{\mu}(a, d)(\mathcal{R}_{\mu}(b, c)(u^{\ast})) + \mathcal{L}^{\ast}_{\mu}(c, a)(\mathcal{M}^{\ast}_{\mu}(b, d)(u^{\ast})) + \mathcal{M}^{\ast}_{\mu}(\mu(b, c, a), d)(u^{\ast})= 0,
\eeq
\beq \label{dtr6}
\mathcal{L}^{\ast}_{\mu}(c, d)(\mathcal{R}^{\ast}_{\mu}(a, b) (u^{\ast})) + \mathcal{R}^{\ast}_{\mu}(a, b)(\mathcal{L}^{\ast}_{\mu}(c, d)(u^{\ast})) + \mathcal{M}^{\ast}_{\mu}(d, a)(\mathcal{M}^{\ast}_{\mu}(c, b)(u^{\ast}))= 0,
\eeq
$ \forall a, b, c, d, x, y, z\in \A, u^{\ast}\in V^{\ast}$.
\end{lemma}
\begin{proposition}
Let $(\mathcal{A}, \mu)$ be a partially associative $3$-ary algebra. Then, $(L_{\mu}, 0, 0, \mathcal{A})$, $(0, 0, R_{\mu}, \mathcal{A})$, $(L_{\mu}, M_{\mu}, R_{\mu}, 
 \mathcal{A})$, $(R^{\ast}_{\mu}, 0, 0, \mathcal{A}^{\ast})$, $(0, 0, 
 L^{\ast}_{\mu}, \mathcal{A}^{\ast})$ and $(R^{\ast}_{\mu}, M^{\ast}_{\mu}, L^{\ast}_{\mu}, 
 \mathcal{A}^{\ast})$are  trimodules of $\mathcal{A}$.
\end{proposition}

\begin{theorem}
 Let $(\mathcal{A}, \mu)$  be a partially associative $3$-ary algebra. Suppose that there is a partially associative $3$-ary algebra structure $\nu$ on its dual space $\mathcal{A}^{\ast}$.

 Then, $(\mathcal{A}, \mathcal{A}^{\ast}, R^{\ast}_{\mu}, M^{\ast}_{\mu}, L^{\ast}_{\mu}, R^{\ast}_{\nu}, M^{\ast}_{\nu}, L^{\ast}_{\nu})$ is a matched pair of partially associative $3$-ary algebras if and only if there is an associative partially $3$-ary  algebra structure on the direct sum $\mathcal{A} \oplus \mathcal{A}^{\ast}$ of the underlying vector spaces of $ \mathcal{A} $ and $ \mathcal{A}^{\ast} $ given by the product $\tau$ defined by
 \beqs
\tau[(x + a^{\ast}), (y + b^{\ast}), (z + c^{\ast})]&=& [\mu(x, y, z) + L^{\ast}_{\nu}(a^{\ast}, b^{\ast})(z) + M^{\ast}_{\nu}(a^{\ast}, c^{\ast})(y) + R^{\ast}_{\nu}(b^{\ast}, c^{\ast})(x)] +\cr
&& [\nu(a^{\ast}, b^{\ast}, c^{\ast}) + L^{\ast}_{\mu}(x, y)(c^{\ast}) + M^{\ast}_{\mu}(x, z)(b^{\ast}) + R^{\ast}_{\mu}(y, z)(a^{\ast})]
\eeqs
for any $x, y, z\in \mathcal{A}, a^{\ast}, b^{\ast}, c^{\ast}\in \mathcal{A}^{\ast}$.
\end{theorem} 
 \section{Associative $3$-ary infinitesimal bialgebras}
We start this section by introducing  the notion of associative $3$-ary infinitesimal bialgebras, which is characterized by a compatible condition between the $3$-ary product and the $3$-ary coproduct. Let us also mention the work by  \cite{[S.Duplij]} on ternary Hopf algebras where a construction of ternary bialgebra is given. 
{
 }
\subsection{Definitions}
\begin{definition}\label{Def5}
A partially associative $3$-ary infinitesimal bialgebra is a triple $(\mathcal{A}, \mu, \Delta)$ such that 
\begin{enumerate}
\item $(\mathcal{A}, \mu)$ is a partially associative $3$-ary algebra,
\item $(\mathcal{A}, \Delta)$ is a partially coassociative $3$-ary coalgebra,
\item $\Delta$ and $\mu$ satisfy the following condition
\beq \textit{\label{cp1}}
\Delta\mu(x, y, z)=(L_{\mu}(x, y)\otimes\id\otimes\id)\Delta(z) + (\id\otimes M_{\mu}(x, z)\otimes\id)\Delta(y) + (\id\otimes\id\otimes R_{\mu}(y, z))\Delta(x).
\eeq
\end{enumerate} 
\end{definition}
\begin{example}\label{ep2}
Let $\mathcal{P}$ be a $2$-dimensionnal vector space with a basis $\lbrace e_{1}, e_{2} \rbrace$. Consider the following products defined by
\beqs
&& \mu: \mathcal{P}\otimes \mathcal{P}\otimes\mathcal{P}\rightarrow \mathcal{P}, \   \mu(e_{1}\otimes e_{1}\otimes e_{1})= e_{2}, \mu(e_{i}\otimes e_{j}\otimes e_{k})=0, i, j, k=1, 2,  (i, j, k)\neq (1, 1, 1); \cr
&& \Delta: \mathcal{P}\rightarrow \mathcal{P}\otimes \mathcal{P}\otimes\mathcal{P}, \  \Delta(e_{1})= e_{2}\otimes e_{2}\otimes e_{2}, \ \Delta(e_{2})= 0.
\eeqs
By a direct computation, we obtain that the triple $(\mathcal{P}, \mu, \Delta)$ is a partially associative $3$-ary infinitesimal bialgebra.
\end{example}

\begin{definition}\label{Def7}
A totally associative $3$-ary infinitesimal bialgebra is a triple $(\mathcal{A}, \mu, \Delta)$ such that 
\begin{enumerate}
\item$(\mathcal{A}, \mu)$ is a totally associative $3$-ary algebra,
\item $(\mathcal{A}, \Delta)$ is a totally coassociative $3$-ary coalgebra,
\item $\Delta$ and $\mu$ satisfy the condition (\ref{cp1}).
\end{enumerate} 
\end{definition}
\begin{example}\label{et2}

Let $\mathcal{T}$ be a $2$-dimensionnal vector space with a basis $\lbrace e_{1}, e_{2}\rbrace$. Consider the following products defined by
\beqs
&& \mu: \mathcal{T}\otimes \mathcal{T}\otimes\mathcal{T}\rightarrow \mathcal{T},\cr
&&\mu(e_{1}\otimes e_{1}\otimes e_{1})= e_{1} \ \ \ \mu(e_{2}\otimes e_{2}\otimes e_{1})= e_{1} + e_{2}\cr
&&\mu(e_{1}\otimes e_{1}\otimes e_{2})= e_{2} \ \ \ \mu(e_{2}\otimes e_{2}\otimes e_{2})= e_{1} + 2e_{2}\cr
&&\mu(e_{1}\otimes e_{2}\otimes e_{1})= e_{2} \ \ \ \mu(e_{1}\otimes e_{2}\otimes e_{2})= e_{1} + e_{2}\cr
&&\mu(e_{2}\otimes e_{1}\otimes e_{1})= e_{2} \ \ \ \mu(e_{2}\otimes e_{1}\otimes e_{2})= e_{1} + e_{2}
\eeqs
and
\beqs
&& \Delta: \mathcal{T}\rightarrow \mathcal{T}\otimes \mathcal{T}\otimes\mathcal{T},\cr
\Delta(e_{1})&=& e_{1}\otimes e_{1}\otimes e_{1} + e_{1}\otimes e_{2}\otimes e_{2} + e_{2}\otimes e_{2}\otimes e_{1} + e_{2}\otimes e_{2}\otimes e_{2} + e_{2}\otimes e_{1}\otimes e_{2}\cr
\Delta(e_{2})&=& e_{1}\otimes e_{1}\otimes e_{2} + e_{1}\otimes e_{2}\otimes e_{2} + e_{2}\otimes e_{1}\otimes e_{1} + e_{2}\otimes e_{2}\otimes e_{1} \cr
&& + e_{2}\otimes e_{1}\otimes e_{2} + e_{1}\otimes e_{2}\otimes e_{1} + 2e_{2}\otimes e_{2}\otimes e_{2}.
\eeqs
By a direct computation, we prove that the triple $(\mathcal{T}, \mu, \Delta)$ is a totally associative $3$-ary infinitesimal bialgebra.
\end{example}
\begin{definition}
A weak totally associative $3$-ary infinitesimal bialgebra is a triple $(\mathcal{A}, \mu, \Delta)$ such that 
\begin{enumerate}
\item$(\mathcal{A}, \mu)$ is a weak totally associative $3$-ary algebra,
\item $(\mathcal{A}, \Delta)$ is a weak totally coassociative $3$-ary coalgebra,
\item $\Delta$ and $\mu$ satisfy the condition (\ref{cp1}).
\end{enumerate} 
\end{definition}

\begin{definition}
Two associative $3$-ary infinitesimal bialgebras $(\mathcal{A}_{1}, \mu_{1}, \Delta_{1})$ and $(\mathcal{A}_{2}, \mu_{2}, \Delta_{2})$ are called equivalent if there exists a vector space isomorphism $f: \mathcal{A}_{1}\rightarrow \mathcal{A}_{2}$ such that
\begin{enumerate}
\item $f: (\mathcal{A}_{1}, \mu_{1})\rightarrow (\mathcal{A}_{2}, \mu_{2})$ is an associative $3$-ary algebra isomorphism, that is, 
\beqs
f\mu_{1}(x, y, z)= \mu_{2}(f(x), f(y), f(z)) \mbox { for all } x, y, z \in \mathcal{A}_{1};
\eeqs 
\item $f: (\mathcal{A}_{1}, \Delta_{1})\rightarrow (\mathcal{A}_{2}, \Delta_{2})$ is a coassociative $3$-ary coalgebra isomorphism, that is, 
\beqs
\Delta_{2}(f(x))= (f\otimes f\otimes f)\Delta_{1}(x) \mbox { for every } x\in \mathcal{A}_{1}.
\eeqs 
\end{enumerate}
\end{definition} 

\begin{example}
Let $\mathcal{P}$ be a $2$-dimensionnal vector space with a basis $\lbrace e_{1}, e_{2}\rbrace$. Consider the following products $\mu_{1}, \mu_{2}: \mathcal{P}\otimes \mathcal{P}\otimes\mathcal{P}\rightarrow \mathcal{P},$\ $\Delta_{1}, \Delta_{2}: \mathcal{P}\rightarrow \mathcal{P}\otimes \mathcal{P}\otimes\mathcal{P},$ defined by
\beqs
&&\mu_{1}(e_{1}\otimes e_{1}\otimes e_{1})= e_{2}, \mu_{1}(e_{i}\otimes e_{j}\otimes e_{k})=0, i, j, k=1, 2,  (i, j, k)\neq (1, 1, 1); \cr
&&\Delta_{1}(e_{1})= e_{2}\otimes e_{2}\otimes e_{2}, \ \Delta_{1}(e_{2})= 0.
\eeqs 
and
\beqs
&&\mu_{2}(e_{2}\otimes e_{2}\otimes e_{2})= e_{1}, \mu_{2}(e_{i}\otimes e_{j}\otimes e_{k})=0, i, j, k=1, 2,  (i, j, k)\neq (2, 2, 2); \cr
&&\Delta_{2}(e_{2})= e_{1}\otimes e_{1}\otimes e_{1}, \ \Delta_{2}(e_{1})= 0.
\eeqs 
By a direct computation, we show that the triples $(\mathcal{P}, \mu_{1}, \Delta_{1})$ and $(\mathcal{P}, \mu_{2}, \Delta_{2})$ are partially associative $3$-ary infinitesimal bialgebras.
 
 Let us consider the linear map $f: \mathcal{P}\rightarrow \mathcal{P}, f(e_{1})= e_{2}, f(e_{2})= e_{1}$. Then, $(\mathcal{P}, \mu_{1}, \Delta_{1})$ and $(\mathcal{P}, \mu_{2}, \Delta_{2})$ are equivalent partially associative $3$-ary infinitesimal bialgebras by the isomorphism $f:(\mathcal{P}, \mu_{1}, \Delta_{1})\rightarrow (\mathcal{P}, \mu_{2}, \Delta_{2})$.
\end{example}
\subsection{Main results}
\begin{theorem}\label{pb1}
Let $(\mathcal{P}, \mu, \Delta)$ be a finite dimensional partially associative $3$-ary infinitesimal bialgebra. Then, 
$(\mathcal{P}^{\ast}, \Delta^{\ast}, \mu^{\ast})$ is a partially associative $3$-ary infinitesimal bialgebra.
\end{theorem}
\textbf{Proof:}
$(\mathcal{P}, \mu, \Delta)$ is a finite dimensional partially associative $3$-ary infinitesimal bialgebra. Then by the Theorem(\ref{teo-co-asso}) and the Theorem(\ref{teo-co-asso1}), $(\mathcal{P}^{\ast}, \Delta^{\ast})$ is a partially associative $3$-ary algebra in the multiplication (\ref{eqco}), and $(\mathcal{P}^{\ast}, \mu^{\ast})$ is a partially coassociative $3$-ary coalgebra in the multiplication (\ref{qdual}). Now, we prove that $\mu^{\ast}: \mathcal{P}^{\ast}\rightarrow \mathcal{P}^{\ast}\otimes \mathcal{P}^{\ast}\otimes \mathcal{P}^{\ast}$ satisfies the identity (\ref{cp1}), that is, the following identity holds for every $\xi, \eta, \gamma \in \mathcal{P}^{\ast} $,
 \beq \textit{\label{cp2}}
\mu^{\ast}(\Delta^{\ast}(\xi, \eta, \gamma))&=&(L_{\Delta^{\ast}}(\xi, \eta)\otimes\id^{\ast}\otimes\id^{\ast})\mu^{\ast}(\gamma) + (\id^{\ast}\otimes M_{\Delta^{\ast}}(\xi, \gamma)\otimes\id^{\ast})\mu^{\ast}(\eta)\cr
&& +  (\id^{\ast}\otimes\id^{\ast}\otimes R_{\Delta^{\ast}}(\eta, \gamma))\mu^{\ast}(\xi)
\eeq
For every $x, y, z\in \mathcal{P}$ and $\xi, \eta, \gamma \in \mathcal{P}^{\ast},$ by identities (\ref{eqco}) and (\ref{qdual}),
\beqs
\langle \mu^{\ast}\Delta^{\ast}(\xi, \eta, \gamma), x\otimes y\otimes z \rangle 
&=&\langle \xi\otimes \eta\otimes \gamma, \Delta\mu(x\otimes y\otimes z)  \rangle \cr
&=&\langle \xi\otimes \eta\otimes \gamma, (L_{\mu}(x, y)\otimes\id\otimes\id)\Delta(z) +\cr
&& (\id\otimes M_{\mu}(x, z)\otimes\id)\Delta(y) +\cr
&& (\id\otimes\id\otimes R_{\mu}(y, z))\Delta(x)\rangle\cr
&=&\langle \xi\otimes \eta\otimes \gamma, (L_{\mu}(x, y)\otimes\id\otimes\id)\Delta(z)\rangle +\cr &&\langle\xi\otimes \eta\otimes \gamma, (\id\otimes M_{\mu}(x, z)\otimes\id)\Delta(y)\rangle +\cr
&&\langle \xi\otimes \eta\otimes \gamma, (\id\otimes\id\otimes R_{\mu}(y, z))\Delta(x)\rangle \cr
&=&\langle \Delta^{\ast}(\id^{\ast}\otimes\id^{\ast}\otimes L^{\ast}_{\mu}(x, y))(\xi\otimes \eta\otimes \gamma), z\rangle +\cr
&& \langle \Delta^{\ast}(\id^{\ast}\otimes M^{\ast}_{\mu}(x, z)\otimes\id^{\ast})(\xi\otimes \eta\otimes \gamma), y\rangle +\cr
&&  \langle\Delta^{\ast}(R^{\ast}_{\mu}(y, z)\otimes\id^{\ast}\otimes\id^{\ast})(\xi\otimes \eta\otimes \gamma), x\rangle\cr
&=&\langle \Delta^{\ast}(\xi\otimes \eta\otimes L^{\ast}_{\mu}(x, y)(\gamma)), z\rangle +\cr
&& \langle \Delta^{\ast}(\xi\otimes M^{\ast}_{\mu}(x, z)(\eta)\otimes \gamma), y\rangle +\cr &&\langle\Delta^{\ast}(R^{\ast}_{\mu}(y, z)(\xi)\otimes \eta\otimes \gamma), x\rangle\cr
&=&\langle L_{\Delta^{\ast}}(\xi, \eta)(L^{\ast}_{\mu}(x, y)(\gamma)), z\rangle +\cr
&& \langle M_{\Delta^{\ast}}(\xi, \gamma)(M^{\ast}_{\mu}(x, z)(\eta)), y\rangle +\cr
&& \langle R_{\Delta^{\ast}}( \eta, \gamma)(R^{\ast}_{\mu}(y, z)(\xi)), x\rangle\cr
&=&\langle L^{\ast}_{\mu}(x, y)(\gamma), L^{\ast}_{\Delta^{\ast}}(\xi, \eta)(z)\rangle +\cr
&& \langle M^{\ast}_{\mu}(x, z)(\eta), M^{\ast}_{\Delta^{\ast}}(\xi, \gamma)(y)\rangle +\cr
&& \langle R^{\ast}_{\mu}(y, z)(\xi), R^{\ast}_{\Delta^{\ast}}(\eta, \gamma)(x)\rangle\cr
&=&\langle \gamma, \mu(x, y, L^{\ast}_{\Delta^{\ast}}(\xi, \eta)(z))\rangle +\cr
&& \langle \eta, \mu(x, M^{\ast}_{\Delta^{\ast}}(\xi, \gamma)(y), z\rangle +\cr
&& \langle \gamma, \mu(R^{\ast}_{\Delta^{\ast}}(\eta, \gamma)(x), y, z)\rangle \cr
&=&\langle \mu^{\ast}(\gamma), x\otimes y\otimes L^{\ast}_{\Delta^{\ast}}(\xi, \eta)(z))\rangle +\cr
&& \langle \mu^{\ast}(\eta), x\otimes M^{\ast}_{\Delta^{\ast}}(\xi, \gamma)(y)\otimes z\rangle +\cr
&& \langle \mu^{\ast}(\gamma), R^{\ast}_{\Delta^{\ast}}(\eta, \gamma)(x)\otimes y\otimes z\rangle\cr
&=&\langle \mu^{\ast}(\gamma), (\id\otimes \id\otimes L^{\ast}_{\Delta^{\ast}}(\xi, \eta))(x\otimes y\otimes z)\rangle +\cr
&& \langle \mu^{\ast}(\eta), (\id\otimes M^{\ast}_{\Delta^{\ast}}(\xi, \gamma)\otimes \id)(x\otimes y\otimes z)\rangle +\cr
&& \langle \mu^{\ast}(\xi), (R^{\ast}_{\Delta^{\ast}}(\eta, \gamma)\otimes \id\otimes \id)(x\otimes y\otimes z)\rangle\cr
&=&\langle (L_{\Delta^{\ast}}(\xi, \eta)\otimes \id^{\ast}\otimes\id^{\ast})\mu^{\ast}(\gamma), x\otimes y\otimes z\rangle + \cr
&&\langle (\id^{\ast}\otimes M_{\Delta^{\ast}}(\xi, \gamma)\otimes \id^{\ast})\mu^{\ast}(\eta), x\otimes y\otimes z\rangle +\cr
&& \langle\id^{\ast}\otimes \id^{\ast}\otimes(R_{\Delta^{\ast}}(\eta, \gamma))\mu^{\ast}(\xi), x\otimes y\otimes z\rangle.
\eeqs
Then, we obtain
\beqs
 \mu^{\ast}\Delta^{\ast}(\xi, \eta, \gamma)&=& (\id^{\ast}\otimes \id^{\ast}\otimes R_{\Delta^{\ast}}(\eta, \gamma))\mu^{\ast}(\xi)+ (\id^{\ast}\otimes M_{\Delta^{\ast}}(\xi, \gamma)\otimes \id^{\ast})\mu^{\ast}(\eta)+\cr
 &&(L_{\Delta^{\ast}}( \xi, \eta)\otimes \id^{\ast}\otimes \id^{\ast})\mu^{\ast}(\gamma).
\eeqs
Hence the identity (\ref{cp2}) holds.

$ \hfill \square $

The triple $(\mathcal{P}^{\ast}, \Delta^{\ast}, \mu^{\ast})$ is called the dual partially associative $3$-ary infinitesimal bialgebra of $(\mathcal{P}, \mu, \Delta)$.

\begin{example}
From Theorem (\ref{pb1}), the dual partially associative $3$-ary infinitesimal bialgebra ($\mathcal{P}^{\ast}, \Delta^{\ast}, \mu^{\ast}$) of ($\mathcal{P}, \mu, \Delta$) in Example (\ref{ep2}) has its multiplications $\Delta^{\ast}: \mathcal{P}^{\ast}\otimes \mathcal{P}^{\ast}\otimes\mathcal{P}^{\ast}\rightarrow \mathcal{P}^{\ast}$ and 

 $ \mu^{\ast}: \mathcal{P}^{\ast}\rightarrow \mathcal{P}^{\ast}\otimes \mathcal{P}^{\ast}\otimes\mathcal{P}^{\ast},$ defined by:
\beqs
&&\Delta^{\ast}(e^{\ast}_{2}\otimes e^{\ast}_{2}\otimes e^{\ast}_{2})= e^{\ast}_{1}, \Delta^{\ast}(e^{\ast}_{i}\otimes e^{\ast}_{j}\otimes e^{\ast}_{k})=0, i, j, k=1, 2, (i, j, k)\neq (2, 2, 2); \cr
&&\mu^{\ast}(e^{\ast}_{2})= e^{\ast}_{1}\otimes e^{\ast}_{1}\otimes e^{\ast}_{1}, \mu^{\ast}(e^{\ast}_{1})=0.
\eeqs 
\end{example}

Similarly, we have the following results.
\begin{theorem}\label{tb1}
Let $(\mathcal{A}, \mu, \Delta)$ be a finite dimensional totally associative $3$-ary infinitesimal bialgebra. Then, 
$(\mathcal{A}^{\ast}, \Delta^{\ast}, \mu^{\ast})$ is a totally associative $3$-ary infinitesimal bialgebra.
\end{theorem}

The triple $(\mathcal{A}^{\ast}, \Delta^{\ast}, \mu^{\ast})$ is   called the dual totally associative $3$-ary infinitesimal bialgebra of $(\mathcal{A}, \mu, \Delta)$.

\begin{example}
From Theorem (\ref{tb1}), the dual totally associative $3$-ary infinitesimal bialgebra ($\mathcal{T}^{\ast}, \Delta^{\ast}, \mu^{\ast}$) of ($\mathcal{T}, \mu, \Delta$) in Example (\ref{et2}) has its multiplications $\Delta^{\ast}: \mathcal{T}^{\ast}\otimes \mathcal{T}^{\ast}\otimes\mathcal{T}^{\ast}\rightarrow \mathcal{T}^{\ast}$ and 

 $ \mu^{\ast}: \mathcal{T}^{\ast}\rightarrow \mathcal{T}^{\ast}\otimes \mathcal{T}^{\ast}\otimes\mathcal{T}^{\ast},$ defined such that $\Delta^{\ast}= \mu$ and $\mu^{\ast}= \Delta$.
\end{example} 
\begin{theorem}
Let $(\mathcal{A}, \mu, \Delta)$ be a  finite dimensional weak totally associative $3$-ary infinitesimal bialgebra. Then, 
$(\mathcal{A}^{\ast}, \Delta^{\ast}, \mu^{\ast})$ is a weak totally associative $3$-ary  bialgebra.
\end{theorem}

The triple $(\mathcal{A}^{\ast}, \Delta^{\ast}, \mu^{\ast})$ is  called the dual weak totally associative $3$-ary infinitesimal bialgebra of $(\mathcal{A}, \mu, \Delta)$.

Now, we study $3$-ary associative bialgebras by means of structure constants.

Let $(\mathcal{A}, \mu, \Delta)$ be a partially $3$-ary associative infinitesimal bialgebra with the multiplications in the basis $e_{1},..., e_{n}$ as follows:
\beq \label{cp3}
\mu(e_{i}, e_{j}, e_{k})=\sum^{n}_{l=1} c^{l}_{ijk}e_{l}, \ \ \Delta(e_{l})=\sum_{1\leq r, s, t\leq n}a^{rst}_{l}e_{r}\otimes e_{s}\otimes e_{t},
\eeq
where $c^{l}_{ijk}, a^{rst}_{l}\in \mathcal{K}, 1\leq i,j, k, l\leq n.$ Then,
\beq \label{cp5}
\Delta\mu(e_{i}, e_{j}, e_{k})=\Delta(\sum^{n}_{l=1} c^{l}_{ijk}e_{l})= \sum^{n}_{l=1} c^{l}_{ijk}\Delta(e_{l})= \sum^{n}_{l=1}\sum_{1\leq r, s, t\leq n}c^{l}_{ijk} a^{rst}_{l}e_{r}\otimes e_{s}\otimes e_{t} 
\eeq
From the identity (\ref{cp1}), we have:
\beq \label{cp4}
\Delta\mu(e_{i}, e_{j}, e_{k})&=&(L_{\mu}(e_{i}, e_{j})\otimes\id\otimes\id)\Delta(e_{k}) + (\id\otimes M_{\mu}(e_{i}, e_{k})\otimes\id)\Delta(e_{j}) +\cr
&& (\id\otimes\id\otimes R_{\mu}(e_{j}, e_{k}))\Delta(e_{i})\cr
&=&(L_{\mu}(e_{i}, e_{j})\otimes\id\otimes\id)\left( \sum_{1\leq r, s, t\leq n}a^{rst}_{k}e_{r}\otimes e_{s}\otimes e_{t}\right)  +\cr
&& (\id\otimes M_{\mu}(e_{i}, e_{k})\otimes\id)\left( \sum_{1\leq r, s, t\leq n}a^{rst}_{j}e_{r}\otimes e_{s}\otimes e_{t}\right)  +\cr
&&  (\id\otimes\id\otimes R_{\mu}(e_{j}, e_{k}))\left( \sum_{1\leq r, s, t\leq n}a^{rst}_{i}e_{r}\otimes e_{s}\otimes e_{t}\right) \cr
&=&\sum_{1\leq r, s, t\leq n}[a^{rst}_{k}\mu(e_{i}, e_{j}, e_{r})\otimes e_{s}\otimes e_{t} + a^{rst}_{j}e_{r}\otimes \mu(e_{i}, e_{s}, e_{k})\otimes e_{t} +\cr
&& a^{rst}_{i}e_{r}\otimes e_{s}\otimes \mu(e_{t},e_{j}, e_{k})] \cr
&=&\sum_{1\leq r, s, t\leq n}\sum_{l=1}^{n}[a^{rst}_{k}c^{l}_{ijr}e_{l}\otimes e_{s}\otimes e_{t} + a^{rst}_{j}c^{l}_{isk}e_{r}\otimes e_{l}\otimes e_{t} + a^{rst}_{i}c^{l}_{tjk} e_{r}\otimes e_{s}\otimes e_{l}]\cr
&=&\sum_{1\leq r, s, t\leq n}\sum_{l=1}^{n}[a^{lst}_{k}c^{r}_{ijl}+ a^{rlt}_{j}c^{s}_{ilk} + a^{rsl}_{i}c^{t}_{ljk}]e_{r}\otimes e_{s}\otimes e_{t}.
\eeq
Comparing the identities (\ref{cp5}) and (\ref{cp4}), we infer
\beq
\sum_{l=1}^{n}[a^{lst}_{k}c^{r}_{ijl}+ a^{rlt}_{j}c^{s}_{ilk} + a^{rsl}_{i}c^{t}_{ljk} -a^{rst}_{l}c^{l}_{ijk}]=0.
\eeq
Conversely, if an  associative $3$-ary algebra $(\mathcal{A}, \mu),$ and a coassociative $3$-ary coalgebra defined by (\ref{cp3}) satisfy the identity (\ref{cp4}), then $\mu$ and $\Delta$ satisfy the identity (\ref{cp1}). Therefore, we have the following result:

\begin{theorem}
Let $\mathcal{A}$ be a vector space with a basis $e_{1},..., e_{n},$ and $(\mathcal{A}, \mu)$ and $(\mathcal{A}, \Delta)$ be associative $3$-ary algebra and coassociative $3$-ary coalgebra, respectively, defined by (\ref{cp3}). Then, $(\mathcal{A}, \mu, \Delta)$ is an associative $3$-ary infinitesimal bialgebra, if and only if $c^{l}_{ijk}$ and $a^{rst}_{l},$ $1\leq i, j, k, l, r, s, t\leq n,$ satisfy identity (\ref{cp4}). 
\end{theorem}

Let $\A$ be an associative 3-ary algebra. Considering  the exchange operator $\sigma: \A\otimes \A\rightarrow \A\otimes\A$ defined as
\beqs
\sigma(a\otimes b)= b\otimes a, \mbox{ for all } a, b\in \A,
\eeqs
 the condition (\ref{cp1}) can be rewritten as:
\beq \label{IATB}
\Delta\mu=(\mu\otimes\id\otimes\id)(\id\otimes\id\otimes\Delta) + (\id\otimes\mu\otimes\id)(\sigma\otimes \id\otimes\sigma)(\id\otimes\Delta\otimes\id)\cr + (\id\otimes\id\otimes\mu)(\Delta\otimes\id\otimes\id).
\eeq
 The relation (\ref{IATB}) can be re-expressed in terms of elements as follows:
\beq\label{inf.}
\Delta\mu(a, b, c)=\sum_{(c)} abc_{1}\otimes c_{2}\otimes c_{3} + \sum_{(b)}b_{1}\otimes ab_{2}c\otimes b_{3} + \sum_{(a)}a_{1}\otimes a_{2}\otimes a_{3}bc, \mbox{ for all } a, b, c\in \A,
\eeq
or, equivalently,
\beq
\Delta\mu(a, b, c)=(L_{\mu}(a, b)\otimes\id\otimes\id)\Delta(c) + (\id\otimes M_{\mu}(a, c)\otimes\id)\Delta(b) + (\id\otimes\id\otimes R_{\mu}(b, c))\Delta(a).
\eeq  
\begin{definition}
An associative 3-ary infinitesimal bialgebra is a triple $(\A, \mu, \Delta)$ consisting of an associative 3-ary algebra $(\A, \mu)$ and a coassociative 3-ary coalgebra $(\A, \Delta)$ such that the condition (\ref{IATB}) holds.
\end{definition}
\begin{proposition}
Let $(\A, \mu, \Delta)$ be an associative 3-ary 
infinitesimal bialgebra. Then so are $(\A, -\mu, \Delta), (\A, \mu, -\Delta), (\A, -\mu, -\Delta)$.
\end{proposition}
\textbf{Proof:}
Using the Eq.(\ref{IATB}), we obtain the results by a direct computation.

$ \hfill \square $

\begin{theorem}
Let $(\A, \mu, \Delta)$ be a  finite dimensional associative 3-ary infinitesimal bialgebra. Then so is $(\A^{\ast}, \Delta^{\ast}, \mu^{\ast})$.
\end{theorem}
\textbf{Proof:}
One has checked that $(\A^{\ast}, \Delta^{\ast})$ is an associative 3-ary algebra  and
$(\A^{\ast} , \mu^{\ast})$ is a coassociative 3-ary coalgebra.  It remains to establish the condition 
(\ref{IATB}) for $\Delta^{\ast}$ and $\mu^{\ast}.$  For that, let us consider $a, b, c\in \A$ and $\varphi, \psi, \chi \in 
\A^{\ast}.$ We  obtain:
\beqs
\langle\mu^{\ast}\circ\Delta^{\ast}(\varphi\otimes\psi\otimes\chi), a\otimes b\otimes c\rangle &=&
\langle \varphi\otimes\psi\otimes\chi, \Delta\circ\mu(a\otimes b\otimes c)\rangle\cr
&=&\langle\varphi\otimes\psi\otimes\chi, (\mu\otimes\id\otimes\id)(\id\otimes\id\otimes\Delta)(a\otimes b\otimes c) \rangle +\cr
&&\langle \varphi\otimes\psi\otimes\chi, (\id\otimes\mu\otimes\id)(\sigma\otimes \id\otimes\sigma)(\id\otimes\Delta\otimes\id)(a\otimes b\otimes c)\rangle\cr
&&+\langle\varphi\otimes\psi\otimes\chi, (\id\otimes\id\otimes\mu)(\Delta\otimes\id\otimes\id)(a\otimes b\otimes c)\rangle\cr
&=&\langle(\id^{\ast}\otimes\id^{\ast}\otimes\Delta^{\ast})(\mu^{\ast}\otimes\id^{\ast}\otimes\id^{\ast})(\varphi\otimes\psi\otimes\chi), a\otimes b\otimes c\rangle +\cr
&&\langle(\id^{\ast}\otimes\Delta^{\ast}\otimes\id^{\ast})(\sigma^{\ast}\otimes\id^{\ast}\otimes\sigma^{\ast})
(\id^{\ast}\otimes\mu^{\ast}\otimes\id^{\ast})(\varphi\otimes\psi\otimes\chi), \cr
&& a\otimes b\otimes c\rangle + \langle(\Delta^{\ast}\otimes\id^{\ast}\otimes\id^{\ast})(\id^{\ast}\otimes\id^{\ast}\otimes\mu^{\ast})(\varphi\otimes\psi\otimes\chi),\cr
&& a\otimes b\otimes c\rangle \mbox{ where } \sigma^{\ast}(\varphi\otimes\psi)= \psi\otimes\varphi,
\eeqs
 showing that the condition (\ref{IATB}) holds for $\Delta^{\ast}$ and $\mu^{\ast}$.
$ \hfill \square $

\section{Concluding remarks}
The main results obtained in this study can be summarized as follows:
\begin{enumerate}
\item[i)] Partially and totally coassociative $3$-ary coalgebras are formulated in Definition \ref{Def1},  Definition \ref{Def2} and Definition \ref{Def3}, and  characterized  in terms of the algebra constant structures in Theorem \ref{theo1}). The relation between a coassociative $3$-ary coalgebra and its dual is given in Theorem \ref{teo-co-asso}).
\item[ii)] Theorem \ref{teo-co-asso1} connects associative $3$-ary algebras to coassociative $3$-ary coalgebras.
\item[iii)] Theorem \ref{theo2} links the isomorphism between two coassociative $3$-ary coalgebras to the isomorphism of their duals.
\item[iv)] The definition of the trimodule of partially associative $3$-ary algebras is given in Definition \ref{Def4}. The matched pairs of totally and partially associative $3$-ary algebras are built in Theorem \ref{Theo. Match. Pair} and Theorem \ref{Theo. Match. Partially}, respectively.
\item[v)] Partially and totally associative $3$-ary infinitesimal bialgebras  are defined in Definition \ref{Def5} and Definition \ref{Def7}, respectively,  while their dual structures are provided in Theorem \ref{pb1} and  Theorem \ref{tb1}. 
\end{enumerate}

 \section*{Aknowledgement}
 This work is partially supported by the Abdus Salam International Centre for Theoretical
 Physics (ICTP, Trieste, Italy) through the Office of External Activities (OEA) - Prj-15. The
 ICMPA is also in partnership with the Daniel Iagolnitzer Foundation (DIF), France.

 \end{document}